\documentclass[preprint,12pt,number,sort&compress]{elsarticle}

\usepackage{graphicx,amssymb,amsthm,amsfonts,subfigure,verbatim}
\usepackage{algorithm,algpseudocode,setspace}
\usepackage{epstopdf,booktabs,color}
\usepackage{makecell}

\usepackage[cmex10]{amsmath}
\usepackage{lineno}

\usepackage{comment}

\def\equationautorefname~#1\null{Eq. (#1)\null}
\def\figureautorefname~#1\null{Fig. #1\null}

\usepackage[textfont=it]{caption}
\usepackage{hyperref}

\journal{Parallel Computing}

\begin{document}

\begin{frontmatter}


\title{A DISTRIBUTED-MEMORY HIERARCHICAL SOLVER FOR GENERAL SPARSE LINEAR SYSTEMS}



\author[sf_icme]{Chao Chen\corref{cor1}}
\ead{cchen10@stanford.edu}
\author[sf_me]{Hadi Pouransari}
\ead{hadip@stanford.edu}
\author[sd]{Sivasankaran Rajamanickam}
\ead{srajama@sandia.gov}
\author[sd]{Erik G. Boman}
\ead{egboman@sandia.gov}
\author[sf_icme,sf_me]{Eric Darve}
\ead{darve@stanford.edu}
\cortext[cor1]{Corresponding author}

\address[sf_icme]{Institute for Computational and Mathematical Engineering, Stanford University, Stanford, CA, USA}
\address[sf_me]{Department of Mechanical Engineering, Stanford University, Stanford, CA, USA}
\address[sd]{Center for Computing Research, Sandia National Laboratories, NM, USA.}

\begin{abstract}
We present a parallel hierarchical solver for general sparse linear systems on distributed-memory machines. For large-scale problems, this fully algebraic algorithm is faster and more memory-efficient than sparse direct solvers because it exploits the low-rank structure of fill-in blocks. Depending on the accuracy of low-rank approximations, the hierarchical solver can be used either as a direct solver or as a preconditioner. The parallel algorithm is based on data decomposition and requires only local communication for updating boundary data on every processor. Moreover, the computation-to-communication ratio of the parallel algorithm is approximately the volume-to-surface-area ratio of the subdomain owned by every processor. We present various numerical results to demonstrate the versatility and scalability of the parallel algorithm.
\end{abstract}


\end{frontmatter}


\section{Introduction}
Solving large sparse linear systems is an important building block -- but often a computational bottleneck -- in many engineering applications. Large sparse linear systems arise, for example, from local discretization of elliptic partial differential equations. Solving a general linear system that has $N$ non-zeros with $O(N)$ computer memory and CPU time is challenging, especially when the underlying physical problem is three-dimensional. Existing methods fall into three categories as follows. 

The first category of methods are sparse direct solvers~\cite{davis2016survey} based on Gaussian elimination. These methods organize computation efficiently and leverage fill-reducing ordering schemes. For example, the nested dissection multifrontal algorithm \cite{george1973nested, duff1986direct} performs elimination according to a specific hierarchical structure of the unknowns. Because of their robustness and efficiency, several state-of-the-art sparse direct solvers have been implemented into software packages, which target sequential \cite{chen2008algorithm, davis2004algorithm, superlu99}, shared-memory \cite{pardiso1, superlumt} and distributed-memory computers \cite{mumps, lidemmel03}. However, sparse direct solvers generally require $O(N^2)$ work and $O(N^{\frac{4}{3}})$ computer memory for a three-dimensional problem of size $N$. This quadratic factorization cost and large memory footprint seriously limit the application of sparse direct solvers to truly large-scale problems. 

The second category consists of iterative solvers. These solvers require only $O(N)$ computer memory to store the linear system and are thus more memory-efficient than sparse direct solvers. They can also achieve the optimal time-complexity if the number of iterations is small. For example, the multigrid method \cite{xu1992iterative} is typically the fastest solver for many discretized elliptic PDEs. However, iterative solvers have three disadvantages. First, the convergence of iterative solvers is not guaranteed for general linear systems. For example, the multigrid method may fail to converge for indefinite linear systems and linear systems coming from coupled PDEs. 
Second, the number of iterations may grow rapidly as the condition number of a linear system increases. Third, the setup phase of some iterative solvers relies on sparse matrix-sparse matrix multiplication, as in the algebraic multigrid method, which is complex to scale~\cite{pplnalu, ipdpsnalu}.

The third category of methods developed for solving sparse linear systems are hierarchical solvers~\cite{grasedyck2009domain, xia2009superfast, ho2015hierarchical, amestoy2015improving, aminfar2016fast,  pouransari2017fast} and their parallel counterparts~\cite{kriemann2005parallel, ghysels2015efficient, wang2016parallel, li2016distributed}. The general idea behind these hierarchical solvers is exploiting the low-rank structure of dense matrix blocks that arise during the elimination process to reduce storage and computational cost. As a result of utilizing the data sparsity, hierarchical solvers, compared with sparse direct solvers, typically have reduced memory footprint and smaller computational complexity. However, some of the existing hierarchical solvers still cannot attain quasilinear complexity for solving three-dimensional problems, and others may be either too complicated to be implemented efficiently or may be restricted to only structured problems.

In this paper, we introduce a new parallel hierarchical solver for solving general sparse linear systems. Our parallel algorithm is based on the sequential algorithm in \cite{pouransari2017fast} called LoRaSp, which computes an approximate factorization with sparse block-triangular factors. In particular, our method eliminates the unknowns cluster by cluster and compresses fill-in blocks in a hierarchical fashion. The singular values of these fill-in blocks are often found to decrease geometrically, and general algebraic techniques, such as the SVD, can be used to compute corresponding low-rank approximations. Since the dropping/truncation rule in our method is based on the decay of singular values, it is expected to be more efficient than other level-based or threshold-based rules, which are typically used in the incomplete LU (ILU) factorization \cite{saad1994ilut}. In our method, the bases computed in low-rank approximations serve the same role as restriction and prolongation operators do in the algebraic multigrid method (AMG) \cite{brandt1986algebraic, stuben2001review}. While the construction of restriction and prolongation operators in AMG may require tuning and manual adjustments for a specific linear system, the low-rank bases in our method are computed in a systematic fashion, regardless of the underlying PDE or physical problem.

There are two differences between our method and other hierarchical solvers. First, while most of hierarchical solvers are developed under either the $\cal H$- \cite{hackbusch1999sparse, hackbusch2000sparse} or the hierarchically semiseparable (HSS) \cite{xia2010fast, chandrasekaran2006fast} matrix frameworks, our method is built upon the ${\cal H}^2$-matrix theory \cite{hackbusch2002data, hackbusch2015mathcal}, which provides a more efficient hierarchical low-rank structure. Under some mild conditions, the computational cost and the memory footprint of our method scale linearly with respect to the problem size, and we observed quasilinear complexity in practice for solving various types of problems. Second, unlike other hierarchical solvers, which are typically combined with the multifrontal algorithm, our method relies on domain partitioning, which naturally leads to a data decomposition scheme in the parallel algorithm. Put another way, these two differences allow our parallel algorithm to have the following three features:
\begin{enumerate}
\item Only local communication is required for every processor.
\item The computation-to-communication ratio is approximately the volume-to-surface-area ratio of the subdomain owned by a processor.
\item The bulk of computation is from using sequential dense linear algebra, which has the potential to be significantly accelerated on modern many-core architectures.
\end{enumerate}

To summarize, this paper presents a parallel hierarchical solver for general sparse linear systems, and especially, our work makes the following three major contributions:
\begin{enumerate}
\item New derivation of the LoRaSp algorithm, which reveals the structure of calculation and data dependency in the original algorithm;
\item Development of a bulk-synchronous parallel algorithm and a task-based asynchronous parallel algorithm with the optimal scheduling strategy;
\item Development of a coloring scheme to extract maximum concurrency in the execution, and discussion on optimizing load-balancing in the presence of coloring constraints.
\end{enumerate}

The remainder of this paper is organized as follows. Section~\ref{sec:serial} presents the sequential algorithm, a new derivation of LoRaSp. Section~\ref{sec:parallel} presents the parallel algorithm, focusing on techniques to keep communication local and to maximize concurrency. Section~\ref{sec:analysis} analyzes the computation and communication cost of the parallel algorithm. Section~\ref{sec:results} presents numerical results to demonstrate the versatility and parallel scalability of our parallel hierarchical solver.

\section{Sequential algorithm} \label{sec:serial}
This section presents our new derivation of the LoRaSp algorithm. Although the algorithm works for a general sparse linear system $Ax = b$, we focus on symmetric positive definite (SPD) systems for ease of presentation. From a high-level perspective, the algorithm computes an approximate factorization of an SPD matrix $A$ with the following steps. 

First, a partitioning of the rows/columns of $A$ is computed algebraically. Suppose $\Pi$ is the set of row/column indices (DOFs). A clustering $\Pi = \cup_i \pi_i$, where $\pi_i$ is a cluster of DOFs, can be computed with a graph partitioner, such as METIS/ParMETIS \cite{karypis1998fast}, Scotch \cite{chevalier2008pt} and Zoltan \cite{ZoltanIsorropiaOverview2012}. Second, some portions of DOFs in every cluster are eliminated as all clusters in $\Pi = \cup_i \pi_i$ are looped over. Specifically, the fill-in blocks associated with a cluster $\pi_s$ are compressed, and $\pi_s$ is split into fine DOFs $\pi_s^f$ and coarse DOFs $\pi_s^c$, where $\pi_s^f$ involves no fill-in. Then $\pi_s^f$ is eliminated safely, which does not propagate any existing fill-in (no level-2 fill-in introduced). Third, after all fine DOFs are eliminated, a smaller linear system, $A_c x_c = b_c$, corresponding to the coarse DOFs remains, and the previous steps are applied repeatedly until $A_c$ is small enough to be factorized with the conventional Cholesky factorization.

\subsection{Low-rank elimination} \label{subsec:sparse}

This subsection illustrates the key step in the algorithm, which exploits the low-rank structure of fill-in blocks to eliminate the fine DOFs in a cluster (if a cluster to be eliminated has no fill-in, e.g., the first eliminated cluster, then all DOFs of this cluster are fine DOFs). Suppose at least one cluster in $\Pi = \cup_i \pi_i$ has been processed (fine DOFs have been eliminated), and consider the Schur complement, $\bar{A}$, which corresponds to the remaining clusters: 
\begin{equation*}
\bar{A} = 
\begin{pmatrix}
A_{ss} & A_{sn} & A_{sw} \\
A_{sn}^T & A_{nn} & A_{nw} \\
A_{sw}^T & A_{nw}^T & A_{ww}
\end{pmatrix}
\end{equation*}
where 
\begin{itemize}
\item $s$ stands for a cluster $\pi_s$ to be processed (``self''),
\item $n$ stands for neighbor clusters (``neighbor''), and
\item $w$ stands for the rest (``well-separated'').
\end{itemize}
Two clusters $\pi_a$ and $\pi_a$ are neighbors if $A(\pi_a,\pi_b) \not = 0$. In other words, the two clusters are connected in the graph of $A$. The definition of \textit{well-separated} is important in ${\cal H}$- and ${\cal H}^2$- matrix theories, which implies a partitioning of the clusters in $n$ and $w$ above. Although our framework is general enough that it can use different definitions of well-separated criteria, here we focus on a simple definition that two clusters are well-separated if they are not neighbors. This definition implies that well-separated blocks, $A_{sw}$ and $A_{ws}=A_{sw}^T$, are fill-in blocks.

Next, fine DOFs in $\pi_s$ will be selected and eliminated. It is observed in \cite{pouransari2017fast} that the fill-in block, $A_{sw}$, is numerically low-rank. With an algebraic method, such as the SVD, we obtain the low-rank approximation
\begin{equation} \label{eqn:lra}
A_{sw} = UZ + O(\epsilon)
\end{equation}
where $U$ is an $m \times k$ orthonormal matrix ($k \ll m$). Using a Gram-Schmidt type of orthogonalization procedure, an orthonormal matrix $V$ of size $m \times (m-k)$ can be computed such that
\begin{equation*}
V^T A_{ss}^{-1} U = 0.
\end{equation*}
We introduce the \textit{sparsification operator}:
\begin{equation} \label{eqn:scaling}
{\cal E}_s =
\begin{pmatrix}
V^T A_{ss}^{-1} &  &  \\
U^T A_{ss}^{-1}& & \\
&I&\\
&&I
\end{pmatrix},
\end{equation}
which decouples $\pi_s^f$ from $\pi_s^c$ and $w$ clusters as
\begin{equation*}
{\cal E}_s \bar{A} {\cal E}_s^T \approx
\begin{pmatrix}
V^TA_{ss}^{-1}V& & V^T A_{ss}^{-1} A_{sn} & \\
& U^TA_{ss}^{-1}U & U^T A_{ss}^{-1} A_{sn} & (U^TA_{ss}^{-1}U)Z\\
A_{sn}^T A_{ss}^{-1} V& A_{sn}^T A_{ss}^{-1} U & A_{nn} & A_{nw}\\
& Z^T (U^TA_{ss}^{-1}U)& A_{nw}^T & A_{ww}
\end{pmatrix}
\end{equation*}
where the first $m-k$ rows/columns correspond to the fine DOFs $\pi_s^f$ and the rest of the DOFs in $\pi_s$ correspond to the coarse DOFs $\pi_s^c$. Note that eliminating $\pi_s^f$ does not propagate any existing fill-in associated with $\pi_s$ in the sense of block elimination and low-rank representation (e.g., $A_{sn}^T A_{ss}^{-1} V$ is not considered as fill-in). In other words, no level-2 fill-in is introduced from the elimination of $\pi_s^f$.

To eliminate $\pi_s^f$, we introduce the \textit{elimination operator}:
\begin{equation} \label{eqn:elim}
{\cal G}_s = 
\begin{pmatrix}
I&&& \\
&I&&\\
- (A_{sn}^T A_{ss}^{-1} V) L^{-T}&&I&\\
&&&I
\end{pmatrix}
\begin{pmatrix}
L^{-1}&&& \\
&I&&\\
&&I&\\
&&&I
\end{pmatrix}
\end{equation}
where $V^T A_{ss}^{-1} V = L L^T$ is the Cholesky factorization, and the first diagonal block in ${\cal G}_s {\cal E}_s \bar{A} {\cal E}_s^T {\cal G}_s^T$ corresponding to $\pi_s^f$ is an identity matrix.

We introduce the \textit{permutation operator}:
\begin{equation}\label{eqn:perm}
P_s =
\begin{pmatrix}
I &&& \\
&&I & \\
&&& I \\
&I&&
\end{pmatrix}
\end{equation}
where the four identity matrices (from the first to the last row) in $P_s$ have the same sizes as the number of DOFs in $\pi_s^f$, $n$ clusters, $w$ clusters and $\pi_s^c$. Note that $P_s$ is used to permute rows and columns corresponding to $\pi_s^c$ to the last in $\bar{A}$.

Combining the above three operators, we define the \textit{low-rank elimination operator}: ${\cal W}_s = P_s {\cal G}_s {\cal E}_s$, and ${\cal W}_s \bar{A} {\cal W}_s^T$ selects and eliminates the fine DOFs in $\pi_s$ (and permutes $\pi_s^c$ to be the last indices). Again, applying ${\cal W}_s$ and ${\cal W}_s^T$ on both sides of $\bar{A}$ does not introduce any level-2 fill-in in $\bar{A}$. The pseudo-code of the low-rank elimination step is shown at lines 15–20 in Algorithm~\ref{alg:lorasp}.

\subsection{Hierarchical solver}
The hierarchical solver repeatedly applies the low-rank elimination procedure to all clusters in $\Pi = \cup_i \pi_i$, and the result is equivalent to computing an approximate factorization of the original matrix, subject to the error of low-rank approximation. The algorithm the hierarchical solver is shown at lines 1–14 in Algorithm~\ref{alg:lorasp}. 

In Algorithm~\ref{alg:lorasp}, $\pi_0$ is the first cluster to be eliminated. Since no fill-in exists yet, ${\cal E}_0$ is an identity matrix, and applying ${\cal W}_0$ and ${\cal W}_0^T$ is one step of standard block Cholesky factorization. This is also true for any cluster that does not involve any fill-in. Suppose the elimination of $\pi_0$ introduces some fill-in for $\pi_1$. Then, the low-rank elimination procedure is applied to $\bar{A}$, the Schur complement in ${\cal W}_0 A {\cal W}_0^T = \begin{pmatrix}I&\\ &\bar{A} \end{pmatrix}$. After applying the low-rank elimination procedure to all clusters in $\Pi$, we are left with the coarse DOFs, $\cup_{i=0}^{m-1} \pi_i^c$, and we can apply the same process to the coarse DOFs.

Algorithm~\ref{alg:lorasp} outputs a factorization of the original matrix $A$, which can be used to solve the corresponding linear system, and the solve phase follows the standard forward and backward substitution, as shown in Algorithm~\ref{alg:lorasp_solve}. 

\begin{algorithm}[hbtp]
  \caption{Hierarchical solver: factorization phase}
  \label{alg:lorasp}   
  \begin{algorithmic}[1]
    \Procedure{Hierarchical\_Factor}{$A$}
    \If{the size of $A$ is small enough}
    \State{Factorize $A$ with the conventional Cholesky factorization}
    \State \Return
    \EndIf
    \State{Partition the graph of $A$ and obtain vertex clusters $\Pi = \cup_{i=0}^{m-1} \pi_i$} 
    \Statex \Comment{{\scriptsize $m$ is chosen to get roughly constant cluster sizes}}
    \State{$A_0 \gets A$}
    \For{$i \gets 0$ \textbf{to} $m-1$}
    \State{$A_{i+1} \gets$ LowRank\_Elimination($A_i$, $\Pi$, $\pi_i$)}
    \EndFor \Comment{{\scriptsize $A_m = {\cal W}_{m-1} \ldots {\cal W}_1 {\cal W}_0 \, A \, {\cal W}_0^T {\cal W}_1^T \ldots {\cal W}_{m-1}^T$}}
    \State{Extract $A_c$ from the block diagonal matrix $A_m \approx \begin{pmatrix}I&\\ & A_c \end{pmatrix}$}
    \Statex \Comment{{\scriptsize $A_c$ is the Schur complement for the coarse DOFs}}
    \State{$A_c^{fac} \gets$ \Call{Hierarchical\_Factor}{$A_c$}} 
    \Statex \Comment{{\scriptsize Recursive call with a smaller matrix}}
    \State \Return 
    $
    A^{fac} =
    {\cal W}_0^{-1} {\cal W}_1^{-1} \ldots {\cal W}_{m-1}^{-1} 
    \begin{pmatrix}
    I & \\
    & A_c^{fac}
    \end{pmatrix}
    {\cal W}_{m-1}^{-T} \ldots {\cal W}_1^{-T} {\cal W}_0^{-T}
    $
    \Statex \Comment{{\scriptsize $A_c^{fac}$ is not written out explicitly}}
    \EndProcedure
    \Statex
    \Procedure{LowRank\_Elimination}{$A_i$, $\Pi$, $\pi_i$}
    \State{Extract $\bar{A}$ from $A_i \approx \begin{pmatrix} I&\\ &\bar{A}\end{pmatrix}$}
    \State{Compute the low-rank elimination operator $\bar{{\cal W}}_i = P_i {\cal G}_i {\cal E}_i$ based on $\bar{A}$}
    \Statex \Comment{{\scriptsize ${\cal E}_i, {\cal G}_i \text{ and } P_i$ are defined in \autoref{eqn:scaling}, \autoref{eqn:elim} and \autoref{eqn:perm}}}
    \State{${\cal W}_i \gets \begin{pmatrix} I&\\ & \bar{{\cal W}}_i \end{pmatrix}$}
    \Statex \Comment{{\scriptsize ${\cal W}_i$ has the same size as $A_i$}}
    \State \Return ${\cal W}_i A_i {\cal W}_i^T$
	\EndProcedure
    \Statex \Comment{{\scriptsize Notation: $a \gets b$ means assign the value $b$ to $a$, whereas $a=b$ means they are equivalent}}
  \end{algorithmic}
\end{algorithm}

\begin{algorithm}[hbtp]
  \caption{Hierarchical solver: solve phase}
  \label{alg:lorasp_solve}   
  \begin{algorithmic}[1]
    \Procedure{Hierarchical\_Solve}{$A^{fac}$, $b$}
    \State{$y \gets$ \Call{Forward\_Substitution}{$A^{fac}$, $b$}}
    \State{$x \gets$ \Call{Backward\_Substitution}{$A^{fac}$, $y$}}
    \State \Return $x$
    \EndProcedure
    \Statex
    \Procedure{Forward\_Substitution}{$A^{fac}$, $b$}
    \State{$y \gets b$}
    \For{$i \gets 0$ \textbf{to} $m-1$}
    \State{$y \gets {\cal W}_i \, y$ \Comment{{\scriptsize $y$ is overwritten}}}
    \EndFor \Comment{{\scriptsize $y = (y_c, y_f)$ is of the concatenation of $y_f$ and $y_c$}}
    \State{Extract $y_f$ and $y_c$ from $y$}
    \Statex \Comment{{\scriptsize $y_f$ and $y_c$ correspond to the fine DOFs and the coarse DOFs}}
    \State{$y_c \gets$ \Call{Forward\_Substitution}{$A_c^{fac}$, $y_c$} \Comment{{\scriptsize $y_c$ is overwritten}}}
    \State \Return $y = (y_f, y_c)$ \Comment{{\scriptsize output the concatenation of $y_f$ and $y_c$}}
    \EndProcedure
    \Statex    
    \Procedure{Backward\_Substitution}{$A^{fac}$, $y$}
    \State{$x \gets y$}
    \For{$i \gets m-1$ \textbf{to} $0$}
    \State{$x \gets {\cal W}_i^T \, x$ \Comment{{\scriptsize $x$ is overwritten}}}
    \EndFor \Comment{{\scriptsize $x = (x_f, x_c)$ is of the concatenation of $x_c$ and $x_f$}}
    \State{Extract $x_f$ and $x_c$ from $x$}
    \Statex \Comment{{\scriptsize $x_f$ and $x_c$ correspond to the fine DOFs and the coarse DOFs}}
    \State{$x_c \gets$ \Call{Backward\_Substitution}{$A_c^{fac}$, $x_c$} \Comment{{\scriptsize $x_c$ is overwritten}}}
    \State \Return $x = (x_f, x_c)$ \Comment{{\scriptsize Output the concatenation of $x_f$ and $x_c$}}
    \EndProcedure
    \Statex \Comment{{\scriptsize Notation: $a \gets b$ means assign the value $b$ to $a$, whereas $a=b$ means they are equivalent}}
  \end{algorithmic}
\end{algorithm}

\paragraph{Relation to LoRaSp}
Algorithm~\ref{alg:lorasp} is mathematically equivalent to LoRaSp \cite{pouransari2017fast}. However, LoRaSp uses \textit{extended sparsification} (introducing auxiliary variables) and needs to maintain a binary tree structure of all DOFs, all of which are avoided in the new derivation. This new derivation also allows more general partitioning strategies. In the special case when recursive bisection is used to compute the partitioning $\Pi = \cup \pi_i$, the new algorithm becomes fairly similar to LoRaSp.

\subsection{Fill-in property} \label{subsec:fillin}
As we emphasized earlier, 
\begin{center}
\textit{level-2 fill-in is never created}
\end{center}
in the hierarchical solver. This fill-in property is not only the fundamental reason why the hierarchical solver is efficient, but also the base of the parallel algorithm. To make the statement of fill-in property precise, we introduce a new notation, $A_B$, to denote the block sparsity pattern of $A$ corresponding to a partitioning of its rows/columns. Given a partition $\Pi=\cup_{i=0}^{m-1} \pi_i$, $A_B$ is an $m$-by-$m$ matrix such that
\begin{equation} \label{eqn:quotient}
A_B(i,j) = 
\begin{cases}
0 & \quad \text{if the subblock $A(\pi_i, \pi_j)=0$,} \\
1 & \quad \text{else}.
\end{cases}
\end{equation}

With the notation of block sparsity pattern, we restate the fill-in property of the hierarchical solver as below:
\begin{center}
\textit{fill-in blocks correspond to non-zeros in $(A_B)^2$ that are not in $A_B$.}
\end{center}
Using ILU terminology, the amount of fill-in never goes beyond ILU(1). Intuitively, in the hierarchical solver, the low-rank elimination operator decouples fine DOFs from fill-in, so eliminating the fine DOFs introduces only level-1 fill-in. See Theorem 5.3 in \cite{pouransari2017fast} for a formal proof.

This fill-in property of the hierarchical solver is the base of the parallel algorithm. Given a partition of all DOFs $\Pi=\cup_i \pi_i$, two clusters $\pi_i$ and $\pi_j$ never ``interact'' if the distance between them is greater than two. The distance is defined as the length of the shortest path between vertex $i$ and vertex $j$ in the graph of $A_B$.

\section{Parallel hierarchical solver} \label{sec:parallel}

Similar to the parallelization of a multigrid method, our parallel algorithm uses a data decomposition scheme, where every processor owns a subdomain of the entire grid. Since every processor exchanges only boundary data with its neighbor processors, the communication is always local in the parallel algorithm. Moreover, the computation-to-communication ratio is approximately the volume-to-surface-area ratio of the subdomain. In other words, the amount of communication is one order of magnitude less than the amount of computation when the subdomain is large enough.

\subsection{Data decomposition}
We present our parallel algorithm for the finest grid, which corresponds to the graph $G$ of $A_B$ (defined in \autoref{eqn:quotient}). The same algorithm applies for the coarse levels. In the following discussion, we use the term ``node'' to mean a cluster of DOFs in $\Pi=\cup_i \pi$, and the term ``edge'' to mean a matrix block (the edge between node $i$ and node $j$ corresponds to block $A(\pi_i, \pi_j)$).

\begin{figure}[hbt]
  \begin{center}
    \scalebox{0.17}{\includegraphics{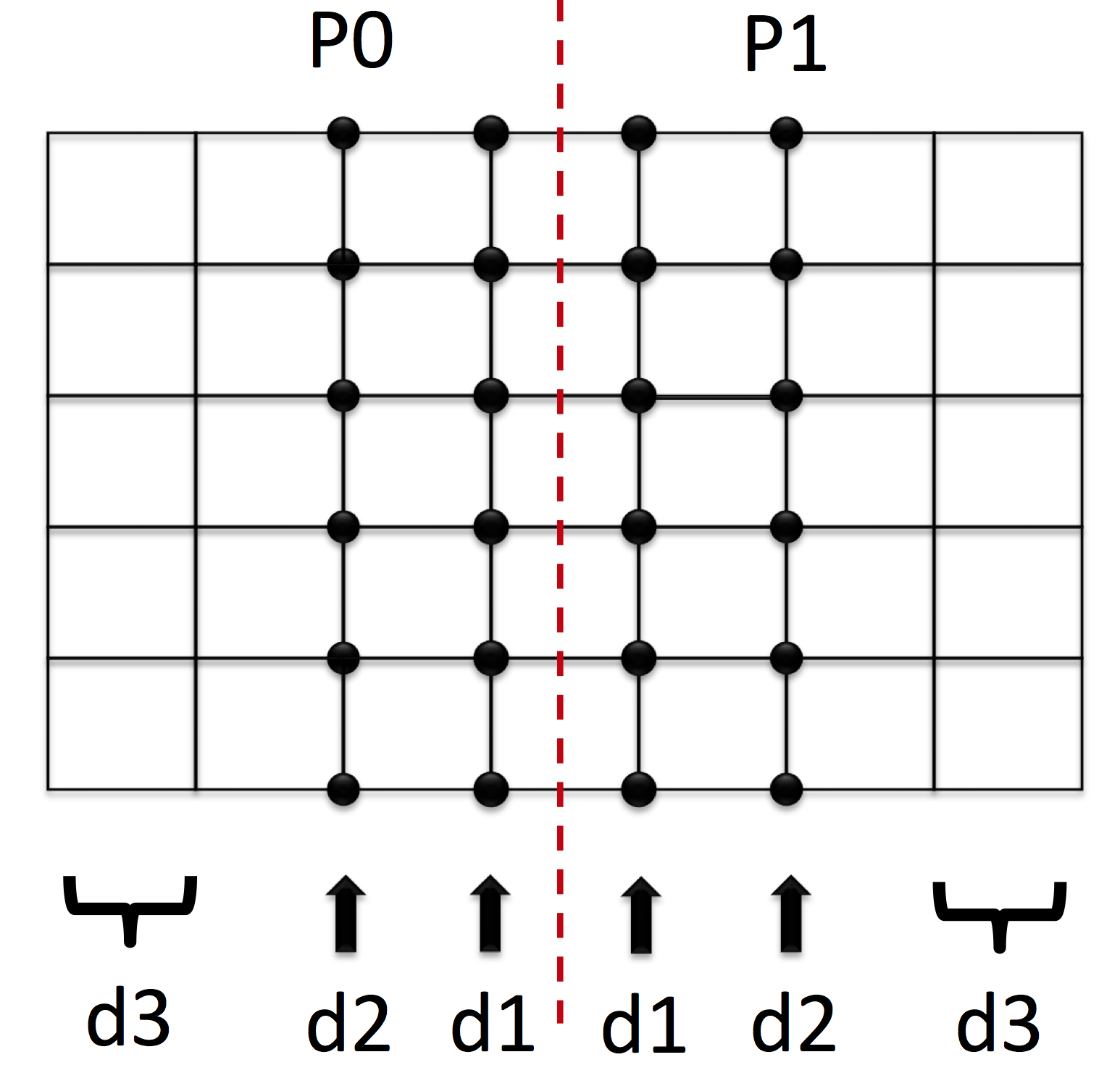}} 
    \caption{A two-processors example showing d1 nodes, d2 nodes and d3 nodes on each processor. The structured grid represents the graph of $A_B$ (defined in \autoref{eqn:quotient}). On each processor, the nodes on the boundary are d1 nodes; d2 nodes are the (local) neighbors of d1 nodes that are not on the boundary; and the remaining interior nodes are d3 nodes.}
    \label{fig:d123}
  \end{center}
\end{figure}

In the parallel algorithm, the global graph $G$ is decomposed among all processors $G = \cup_P G_P$, where processor $P$ owns a subgraph $G_P$. With this decomposition, nodes owned by one processor can be classified into three categories:
\begin{itemize}
\item d1 nodes: boundary nodes,
\item d2 nodes: (local) neighbors of d1 nodes that are not on the boundary,
\item d3 nodes: the remaining interior nodes.
\end{itemize}
Fig. \ref{fig:d123} shows a simple example to clarify the above definitions. We assume that the matrix $A$ is distributed by rows among all processors. Every processor $P$ owns the submatrix corresponding to the local graph $G_P$ and also stores the edges to neighbor processors. For example, the matrix $A_0$ owned by processor $P0$ in Fig. \ref{fig:d123} is
\[
A_0  = 
\begin{pmatrix}
A_0^{d3} & A_0^{d3,d2} && \\
A_0^{d2,d3}  & A_0^{d2} & A_0^{d2,d1} &\\
 &  A_0^{d1,d2} & A_0^{d1} & A_{0,1}^{d1}
\end{pmatrix}
\]
where for the sake of simplicity we assumed that d3 nodes were ordered before d2 nodes and d1 nodes. The coupling with processor $P1$ is contained within the $A_{0,1}^{d1}$ block, since only d1 nodes have edges to processor $P1$.

\subsection{Bulk-synchronous parallel algorithm}

Recall the discussion in Section~\ref{subsec:fillin} that no level-2 fill-in exists in the algorithm. The implication is that we can process (apply low-rank elimination operators on) two clusters/nodes $\pi_i$ and $\pi_j$ in parallel, as long as the distance between them is greater than two. The distance is defined as the length of the shortest path between vertex $i$ and vertex $j$ in the graph of $A_B$. 

This fill-in property motivates us to design the parallel algorithm as follows. First, since the distance between a local d3 node and any remote node is at least three, processing d3 nodes is embarrassingly parallel.

Second, since the distance between a local d2 node and a remote d2 node is at least three, d2 nodes can be processed in parallel. However, communication is required because the shortest distance between a local d2 node and a remote d1 node is two, and the amount of communication is proportional to the number of d2 nodes.

Last, since d1 nodes are coupled between neighbor processors, we color them in a way such that d1 nodes with the same color can be processed in parallel; the coloring scheme is discussed in Section \ref{subsec:coloring}. We introduce two notations: ${\cal N}_1(P)$ for the neighbors of processor $P$ and ${\cal N}_2(P) = {\cal N}_1({\cal N}_1(P)) \backslash \{P\}$, i.e., neighbors of the neighbors of processor $P$. Since the distance between a d1 node on processor $P$ and a remote d1 node owned by a processor in ${\cal N}_2(P)$ can be two, they may need to communicate. The total amount of communication is proportional to the number of d1 nodes. 

To be more specific, the communication in the parallel algorithm corresponds to the ``right-looking style'' communication in a conventional sparse direct solver, and two types of communication are needed:
\begin{enumerate}
\item Suppose we eliminate the fine DOFs of a d1 node on processor $P$. The resulting Schur complement may contain (1) edges that connect a local node on $P$ and a remote node on $P$'s neighbor, and (2) edges that connect two remote d1 nodes on two different processors in ${\cal N}_2(P)$. Therefore, every processor communicates only with its neighbors.
\item Consider the application of the sparsification operator (\autoref{eqn:scaling}). A sparsification of a (local) d2 node may lead to communication with a neighbor processor, and
a sparsification of a (local) d1 node may involve communication with a processor in ${\cal N}_2(P)$.
\end{enumerate}

The computation-to-communication ratio is closely related to the geometry of the subdomain on every processor. The communication volume is proportional to the number of d1 nodes (note in particular that local d2 nodes only exchange data with remote d1 nodes, never with remote d2 nodes), while the amount of computation is proportional to the total number of d1 nodes, d2 nodes and d3 nodes. Therefore, the computation-to-communication ratio is approximately the volume-to-surface-area ratio of the subdomain owned by every processor.

The above discussion is summarized in Algorithm \ref{alg:par1}\footnote{The pseudo-code is in \textit{single program multiple data} (SPMD) style.} (some details are skipped, and the focus is on the high level idea of processing d1 nodes, d2 nodes and d3 nodes separately).

\begin{algorithm}[htbp]
  \caption{Bulk-synchronous parallel algorithm}
  \label{alg:par1}
  \begin{algorithmic}[1] 
    \small
  	\Procedure{HSolver\_Bulk\_Synchronous}{local vertex clusters $\Pi$, local subgraph $G$, local submatrix $A$}
    \Statex \Comment{{\scriptsize synchronous communication is used in this algorithm}}
    \Statex \Comment{{\scriptsize this is processor $P$}}
    \State{Partition $\Pi$ into d1 nodes $\Pi^{(1)}$, d2 nodes $\Pi^{(2)}$, and d3 nodes $\Pi^{(3)}$}
    \State{Parallel (distance-2) coloring of d1 nodes} \Comment{{\scriptsize discussed in Section~\ref{subsec:coloring}}}
    \Statex
    \For{color $i$ = 1 \textbf{to} NumColors} \Comment{{\scriptsize start processing d1 nodes}} 
    \For{$\pi_j \in \Pi^{(1)}$ with color $i$}
    	\State{$A \gets$ LRE($A$, $\Pi$, $\pi_j$)}
    \EndFor 
    \State{Communicate with processors in ${\cal N}_2(P)$}
    \Comment{{\scriptsize neighbor of neighbor processors}}
    \EndFor \Comment{{\scriptsize finish processing d1 nodes}}
    \Statex 
    \For{$\pi_i \in \Pi^{(2)}$}  \Comment{{\scriptsize start processing d2 nodes}}
    \State{$A \gets$ LRE($A$, $\Pi$, $\pi_i$)}
    \EndFor
	\State{Communicate with processors in ${\cal N}_1(P)$} \Comment{{\scriptsize neighbor processors}}
    \Statex \Comment{{\scriptsize finish processing d2 nodes}}
    \For{$\pi_i \in \Pi^{(3)}$} \Comment{{\scriptsize start processing d3 nodes}}
    \State{$A \gets$ LRE($A$, $\Pi$, $\pi_i$)}
    \EndFor \Comment{{\scriptsize finish processing d3 nodes}}
    \Statex \Comment{{\scriptsize no communication needed for d3 nodes}}
	\State{Extract $\Pi^C$, $G^C$, $A^C$ associated with the coarse level}
	\State \Call{HSolver\_Bulk\_Synchronous}{$\Pi^C$, $G^C$, $A^C$}
	\EndProcedure
	\singlespacing
  \end{algorithmic}
\end{algorithm}

\subsection{Task-based asynchronous parallel algorithm}

Algorithm \ref{alg:par1} uses a classical bulk-synchronous style of programming. In this subsection, we explore an asynchronous ``task-based'' approach. Consider the critical path of execution in Algorithm \ref{alg:par1}: d1 nodes must be processed first, and d2 nodes should be processed next. Therefore, the asynchronous algorithm always attempts to process d1 nodes; if none is ready, d2 nodes are processed; and finally, if all else fail, d3 nodes are processed. Processing d3 nodes does not require communication, which is the fall-back if no d1 node or d2 node can be processed. The pattern of communication is organized as follows.

\medskip
\noindent Attempt to process a d1 node $\pi^{(1)}$ with color $i$ on processor $P$ as follows.
\begin{enumerate}
\item Check that the communication for $\pi^{(1)}$ is complete, which is associated with remote d1 nodes owned by ${\cal N}_2(P)$ that have colors smaller than $i$.
\item Apply the sparsification operator (\autoref{eqn:scaling}) on $\pi^{(1)}$.
\item {\bf Send} the updated edges to remote d1 nodes owned by ${\cal N}_2(P)$, and to remote d2 nodes owned by ${\cal N}_1(P)$.
\item Eliminate the fine DOFs in $\pi^{(1)}$.
\item {\bf Send} the updated edges (in the Schur complement) to remote d1 nodes owned by ${\cal N}_1(P)$.
\end{enumerate}

\medskip
\noindent Otherwise, attempt to process a d2 node $\pi^{(2)}$ as follows.
\begin{enumerate}
\item Check that the communication for $\pi^{(2)}$ is complete, which is associated with remote d1 nodes owned by ${\cal N}_1(P)$.
\item Apply the sparsification operator (\autoref{eqn:scaling}) on $\pi^{(2)}$.
\item {\bf Send} the updated edges to remote d1 nodes owned by ${\cal N}_1(P)$ (this data is only used at the coarse level).
\item Eliminate the fine DOFs in $\pi^{(2)}$.
\end{enumerate}

\medskip
\noindent Otherwise, process a d3 node.

The above discussion is summarized in Algorithm~\ref{alg:async_par}\footnote{The pseudo-code is in SPMD style.}, which is based on a critical path analysis of the execution. This algorithm executes the phase with highest priority first (located along the critical path), followed by a phase with medium priority, and low priority. Communication is initiated in an asynchronous manner as soon as the data is available to achieve the maximum concurrency. The large number of executable tasks allows hiding some of the communication latency, thereby minimizing the idle time due to communication. 
 
\begin{algorithm}[htbp]
  \caption{Asynchronous parallel algorithm}
  \label{alg:async_par}  
  \begin{algorithmic}[1]
    \small \onehalfspacing    
  	\Procedure{HSolver\_Asynchronous}{local vertex clusters $\Pi$, local subgraph $G$, local submatrix $A$}
	\Statex \Comment{{\scriptsize asynchronous communication is used in this algorithm}}
    \Statex \Comment{{\scriptsize this is processor $P$}}
    \State{Partition $\Pi$ into d1 nodes $\Pi^{(1)}$, d2 nodes $\Pi^{(2)}$, and d3 nodes $\Pi^{(3)}$}
    \State{Parallel (distance-2) coloring of d1 nodes} \Comment{{\scriptsize discussed in Section~\ref{subsec:coloring}}}
	\State{color $i \gets 1$} \Comment{{\scriptsize start with the smallest color}}
	\While{$\Pi^{(1)} \cup \Pi^{(2)} \cup \Pi^{(3)}$ is not empty}
    \Statex \Comment{{\scriptsize attempt d1 nodes first}}
    \If{$\pi_j \in \Pi^{(1)}$ has color $i$ and communication for $\pi_j$ is complete}
    	\State{$A \gets$ LRE($A$, $\Pi$, $\pi_j$)}
		\State{Send data to ${\cal N}_2(P)$ processors}
	    \Comment{{\scriptsize for remote nodes with colors $> i$}}
	    \State{Pop $\pi_j$ from $\Pi^{(1)}$}
		\Statex \Comment{{\scriptsize attempt d2 nodes secondly}}
	\ElsIf{communication for $\pi_j \in \Pi^{(2)}$ is complete}
	    \State{$A \gets$ LRE($A$, $\Pi$, $\pi_j$)}
		\State{Send data to ${\cal N}_1(P)$ processors}
        \Comment{{\scriptsize for remote d1 nodes}}
	    \State{Pop $\pi_j$ from $\Pi^{(2)}$}
		\Statex \Comment{{\scriptsize attempt d3 node last}}
	\ElsIf{$\pi_j \in \Pi^{(3)}$}
		    \State{$A \gets$ LRE($A$, $\Pi$, $\pi_j$)}
    	    \State{Pop $\pi_j$ from $\Pi^{(3)}$}
    \EndIf
    \State{If all d1 nodes with color $i$ have been processed, $i \gets i+1$}
	\Statex \Comment{{\scriptsize processing d1 nodes with color $i$ is on the critical path of execution}}
    \EndWhile
	\State{Wait for remaining communication}
	\State{Extract $\Pi^C$, $G^C$, $A^C$ associated with the coarse level}
	\State \Call{HSolver\_Asynchronous}{$\Pi^C$, $G^C$, $A^C$}
	\EndProcedure
	\singlespacing	
  \end{algorithmic}
\end{algorithm}

\subsection{Coloring of d1 nodes} \label{subsec:coloring}

The d1 nodes are coupled between neighbor processors. To maximize concurrency, we use a graph coloring scheme to assign colors to all d1 nodes such that nodes with the same color can be processed in parallel. 
For that purpose, a pair of d1 nodes $\pi_i$ and $\pi_j$ must have different colors, if the following two conditions hold. 
\begin{itemize} \label{eqn:d1_color}
\item $\pi_i$ and $\pi_j$ are owned by different processors.
\item the distance between $\pi_i$ and $\pi_j$ is less than or equal to two. 
\end{itemize}

\begin{figure}[htbp]
  \begin{center}
    \scalebox{0.25}{\includegraphics{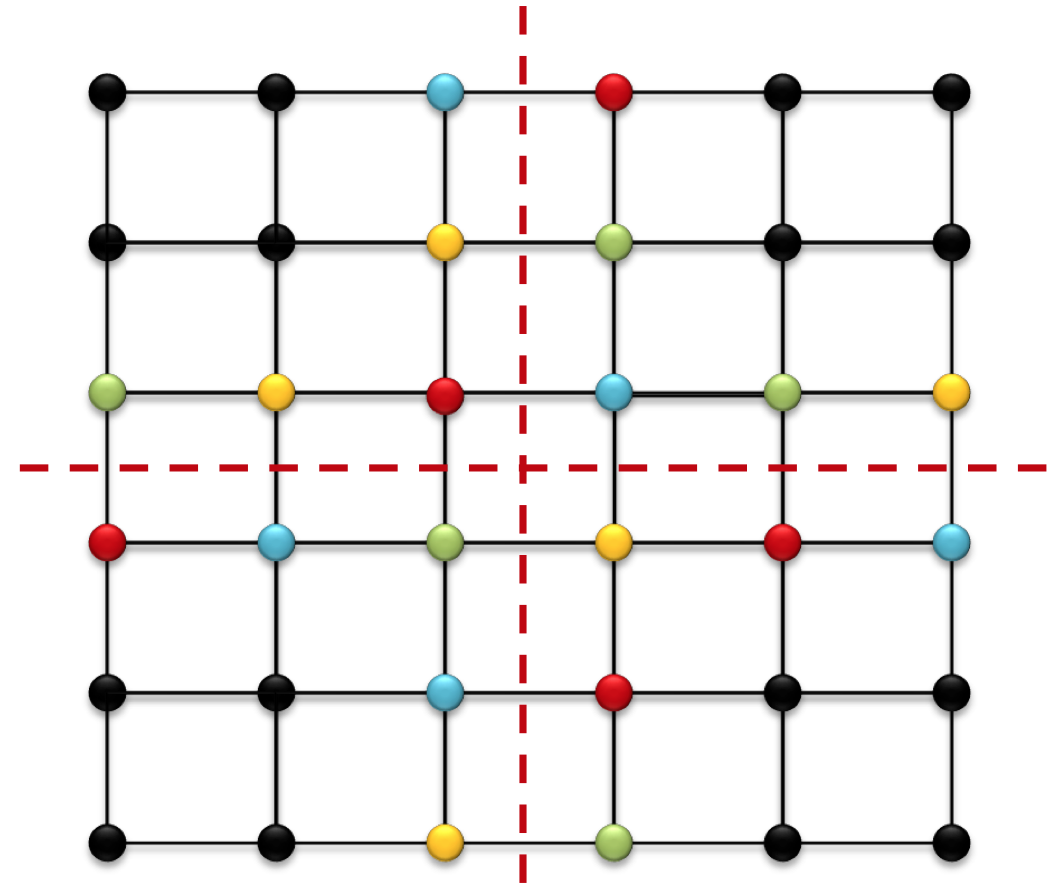}} 
    \caption{A four-processors example of the coloring, where the grid is split onto four processors. A pair of nodes that are at a distance less than two from each other must have different colors, unless they are owned by the same processor. With this coloring result, all four processors are able to process a subset of the d1 nodes concurrently. As seen in this example, however, load imbalance is hard to reduce. Each processor has one node for three given colors and two nodes for the last fourth color, leading to a factor of 2/1 load imbalance during the loop over colors.}
    \label{fig:coloring}
  \end{center}
\end{figure}

For best performance such a coloring should both minimize the number of colors and be load-balanced (the number of nodes of a given color should be roughly constant across all processors). However, even with only the first objective, the coloring problem is NP-hard. Fortunately, fast linear-time greedy heuristics work well in practice for graph coloring. One option is to bring the boundary graph to a single processor and compute the coloring sequentially, but this is not scalable in terms of memory. We therefore use the parallel distance-2 coloring~\cite{bozdag2010distributed} routine in Zoltan~\cite{ZoltanIsorropiaOverview2012}. One example of the coloring result is shown in \autoref{fig:coloring}.

\section{Complexity analysis}\label{sec:analysis}
In this section, we first review the computational cost and the memory consumption of LoRaSp, and clarify the assumptions for its linear complexity. Based on these results, we then analyze the parallel algorithm in terms of the computation complexity, communication complexity and memory consumption. The computational cost and the memory consumption of LoRaSp have been analyzed in \cite{pouransari2017fast}. Below, we rephrase Theorem 5.4 in \cite{pouransari2017fast}, which summarizes the results.

\paragraph{Linear complexity conditions} The computational cost of the factorization and the memory consumption in LoRaSp scales as $O(N r^2)$ and $O(N r)$, respectively, with respect to the problem size $N$ if the following two conditions hold. First, for every cluster of DOFs, the number of neighbors is bounded by a constant. Second, the largest cluster size at the first level (level 0), $r$, is bounded by a constant; and $r_i$, the largest cluster size at level $i$, satisfies the relationship that $r_i < \alpha^i \, r \,(i>0)$, where $0 < \alpha < 2^{1/3}$.

For the rest of this section, we assume the linear complexity conditions always hold. We also assume the linear system is evenly distributed among all processors in the parallel algorithm. Since the computation and the storage are evenly split across all processors, the computational cost and memory consumption are $O(N r^2 / p)$ and $O(N r / p)$ on every processor, where $p$ is the number of processors. 

In the parallel algorithm, communication happens at two places. First, all processors exchange boundary data until the problem size becomes small enough to be factorized efficiently with the sequential Cholesky factorization. Second, the small matrix is gathered to a single processor with a reduction tree and the sequential Cholesky factorization is performed.

For every processor, the amount of communication is dominated by that from the first phase, which is proportional to the number of boundary clusters. Since the matrix associated with a cluster has at most $r^2$ entries, the total amount of communication is 
\[
O\bigg( r^2 \cdot \big(\frac{N}{rp} \big)^{2/3} \bigg) 
= O\bigg( \big(\frac{N r^2}{p} \big)^{2/3} \bigg)
\]
where $\frac{N}{rp}$ is the number of clusters on every processor, and $(\frac{N}{rp})^{2/3}$ is the surface area (assuming an underlying 3D subdomain).

The number of messages sent by every processor in the first communication phase is roughly the number of colors at all levels of grids. If the number of colors at every level is bounded by a constant, then the number of messages is proportional to the number of levels, i.e., $O(\log(N/(rp)))$. In the second communication phase, the number of messages needed is $O(\log p)$ for the conventional Cholesky factorization. Therefore, the total number of messages sent by every processor is 
\[
O\bigg(\log\big(\frac{N}{rp}\big)\bigg) + O(\log p).
\]

To measure parallel scalability, we use standard definitions of parallel speedup, strong scaling efficiency and weak scaling efficiency as follows. Let $T(N, p)$ be the wall-clock time to factor a matrix of size $N$ on $p$ processors. The parallel speedup is defined as 
\begin{equation} \label{eqn:speedup}
S(N,p) = \frac{T(N, p_0)}{T(N, p)}
\end{equation}
where $p_0$ is the smallest number of processors on which the baseline is obtained. The efficiency of strong scaling (the total problem size is fixed as the number of processors increases) is 
\begin{equation} \label{eqn:strong}
E_s = \frac{S(N, p)\times p0}{p}.
\end{equation}
The efficiency of weak scaling (the total problem size increases proportionally to the number of processors, or the problem size per processor is fixed) is
\begin{equation} \label{eqn:weak}
E_w = \frac{T(N\,p_0, p_0)}{T(N \, p, p)}.
\end{equation}

Following the above definitions, we can derive corresponding formulas for our parallel algorithm. As an example, the weak scaling efficiency of our parallel algorithm is
\[
E_w = \frac{N}{N+\log N + \log p}
\]
where $N$ is the problem size per processor, $p_0=1$ is assumed, and the cluster size $r$, a constant, is not shown.
Therefore, a constant efficiency $E_w$ with a fixed problem size per processor should not be expected, even if the parallel algorithm is implemented perfectly with perfect underlying hardware (no network saturation, fixed diameter, etc.). In our case, maintaining a fixed $E_w$ implies that the total problem size should increase at least as fast as $O(p \log(p))$\footnote{$O(p \log(p))$ is the iso-efficiency function of our parallel algorithm~\cite{Grama1993}. This efficiency result is similar to that of computing parallel reduction.}.

\section{Numerical results} \label{sec:results}

This section presents various results to show the versatility and parallel scalability of the asynchronous parallel algorithm in Algorithm~\ref{alg:async_par}. Our implementation is based on the algorithm in Section~\ref{sec:serial} and the LoRaSp algorithm. To run on distributed-memory machines, our code is written with the {\it message passing interface} (MPI). Sequential results are obtained by running our code with one MPI rank. Unless otherwise stated, all tests were run at the NERSC Edison supercomputer\footnote{http://www.nersc.gov/users/computational-systems/edison/configuration/}. Each node of Edison has two 12-core Intel ``Ivy Bridge" processors and nodes are connected by a Cray Aries high-speed interconnect with Dragonfly topology.

\subsection{Linear scalability results} \label{subsec:results_serial}

In this subsection, we present results for solving sequences of problems to demonstrate the scalability of the hierarchical solver. In practice, the linear complexity conditions in \autoref{sec:analysis} may not hold, but we shall see that the total cost of factorization and solve scales almost linearly as the problem size increases.

The results presented in this subsection are the number of iterations and the total time (factorization + solve) for solving three PDEs (discretized with the seven-point stencil on the unit cube). The first is Poisson's equation:
\[
-\Delta u(x) = f(x).
\]

The second is variable-coefficient Poisson's equation:
\[
-\nabla \cdot ( a(x) \nabla u(x) ) = f(x),
\]
where $a(x)$ is a quantized high-contrast random field generated in the following way: (1) initialize $a(x)$ randomly with a uniform distribution; (2) convolve the initial $a(x)$ with Gaussian distribution of deviation $4h$, where $h$ is the stencil spacing; and (3) set $a(x)$ to $10^2$ if it is larger than $0.5$, otherwise $10^{-2}$. We chose a quantized high-contrast random field because these are problems known to be difficult to solve using iterative methods. Although the performance worsens with our hierarchical solver, the algorithm still remains very efficient. For Poisson's equation and variable-coefficient Poisson's equation, we use our solver as a preconditioner for CG with a tolerance of $10^{-12}$.

The third is the Helmholtz equation:
\[
(-\Delta - k^2) u(x) = f(x),
\]
where $k$ is the wave number. We fix the number of DOFs per wave length and increased $k$ proportionally to the number of DOFs in each dimension. As the frequency increases, hierarchical solvers will eventually break down because the ranks required to reach a good accuracy become too large for the method to be efficient. We fixed the resolution to 32 DOFs per wavelength, and the frequencies increased from $f=1$Hz ($32^3$ grid) to $f=4$Hz ($128^3$ grid). For the Helmholtz equation, we use our solver as a preconditioner for GMRES with a lower tolerance of $10^{-3}$ because the Helmholtz equation is often solved with relatively low accuracy in applications such as seismic imaging. 

\begin{figure}[htbp]
  \centering
  \subfigure[Iteration number (Poisson)]
  {\includegraphics[width=62mm]{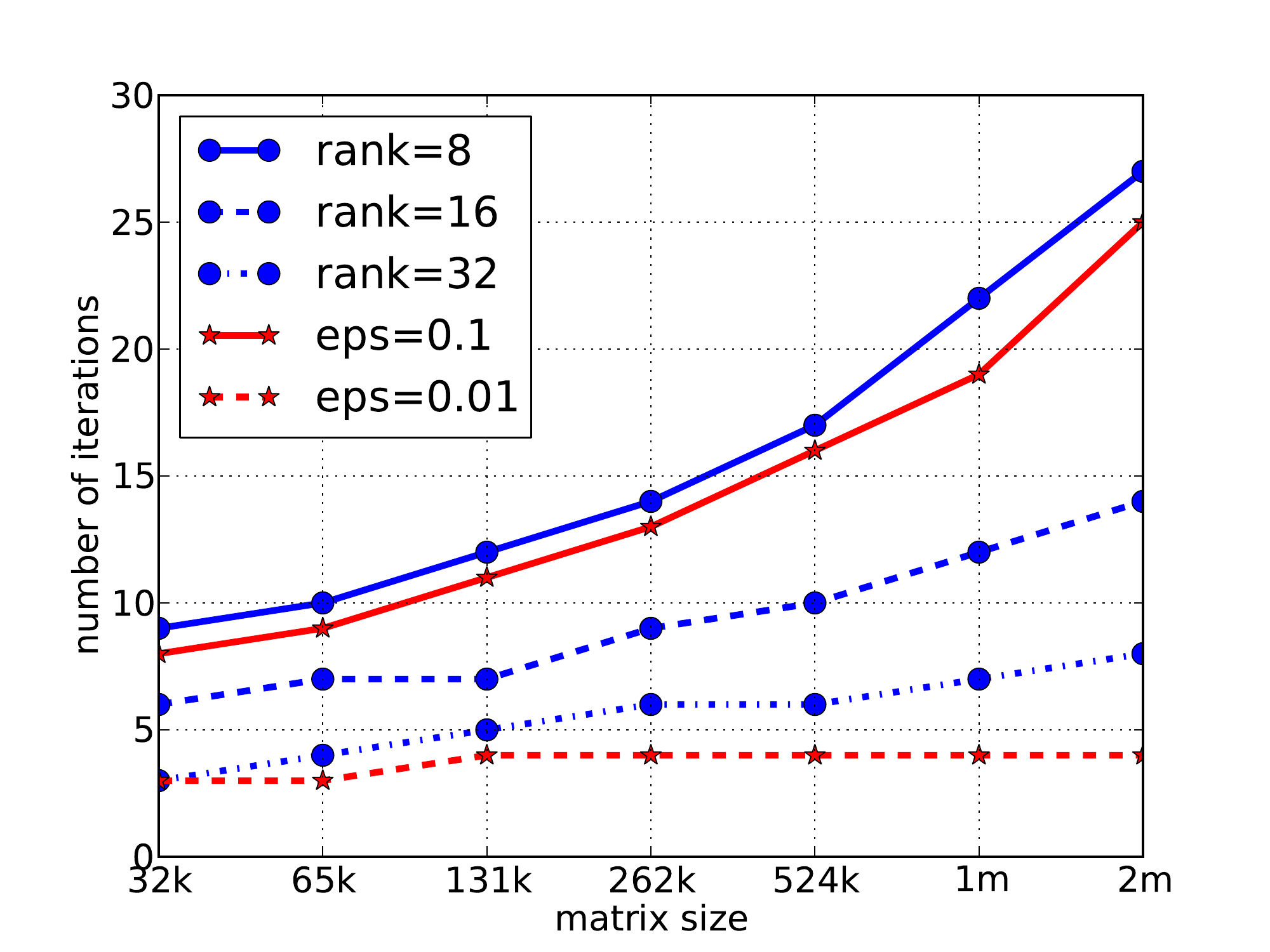}}
  \subfigure[Total time (Poisson)]
  {\includegraphics[width=62mm]{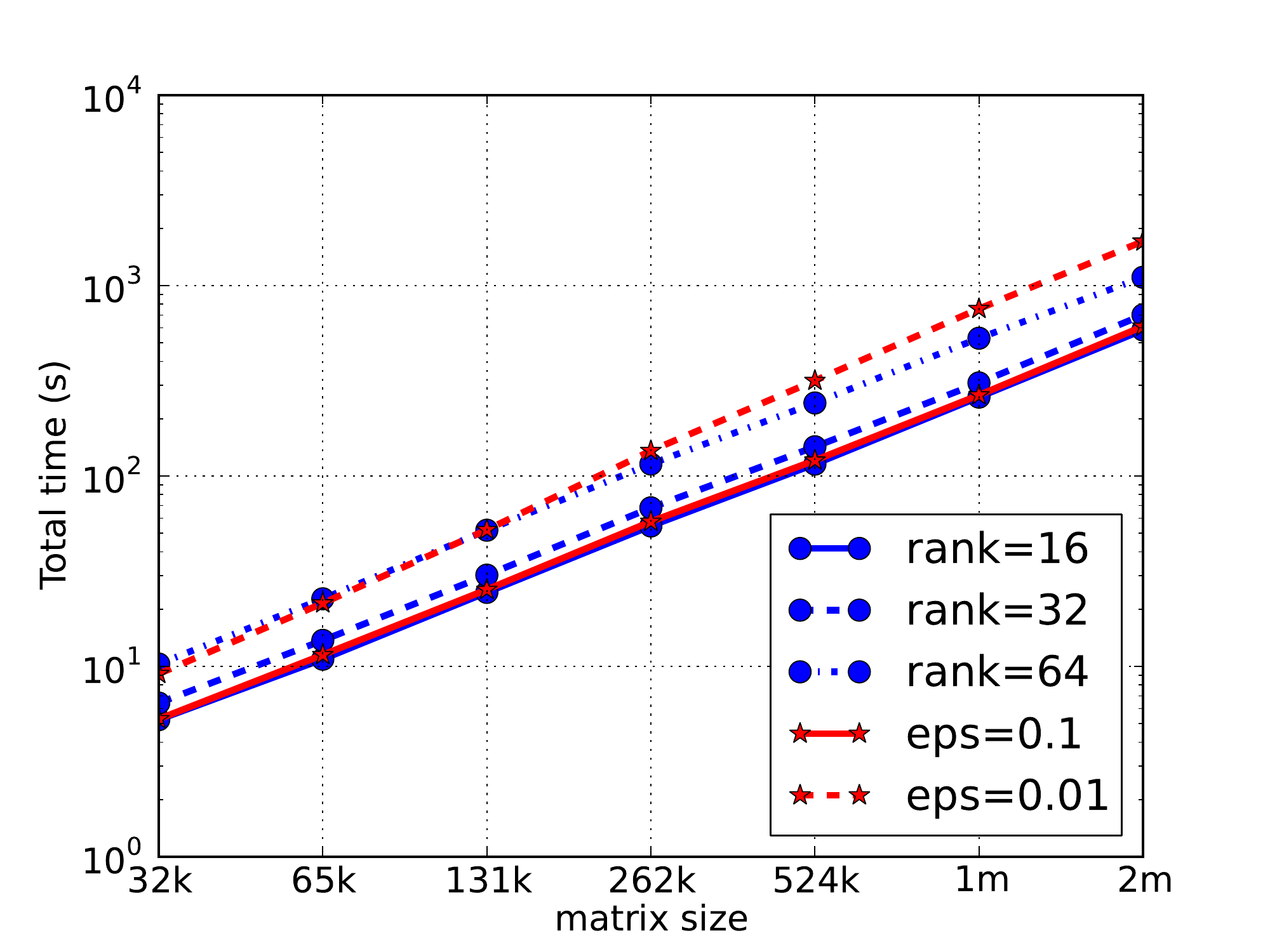}}
  \subfigure[Iteration number (VC-Poisson)]
  {\includegraphics[width=62mm]{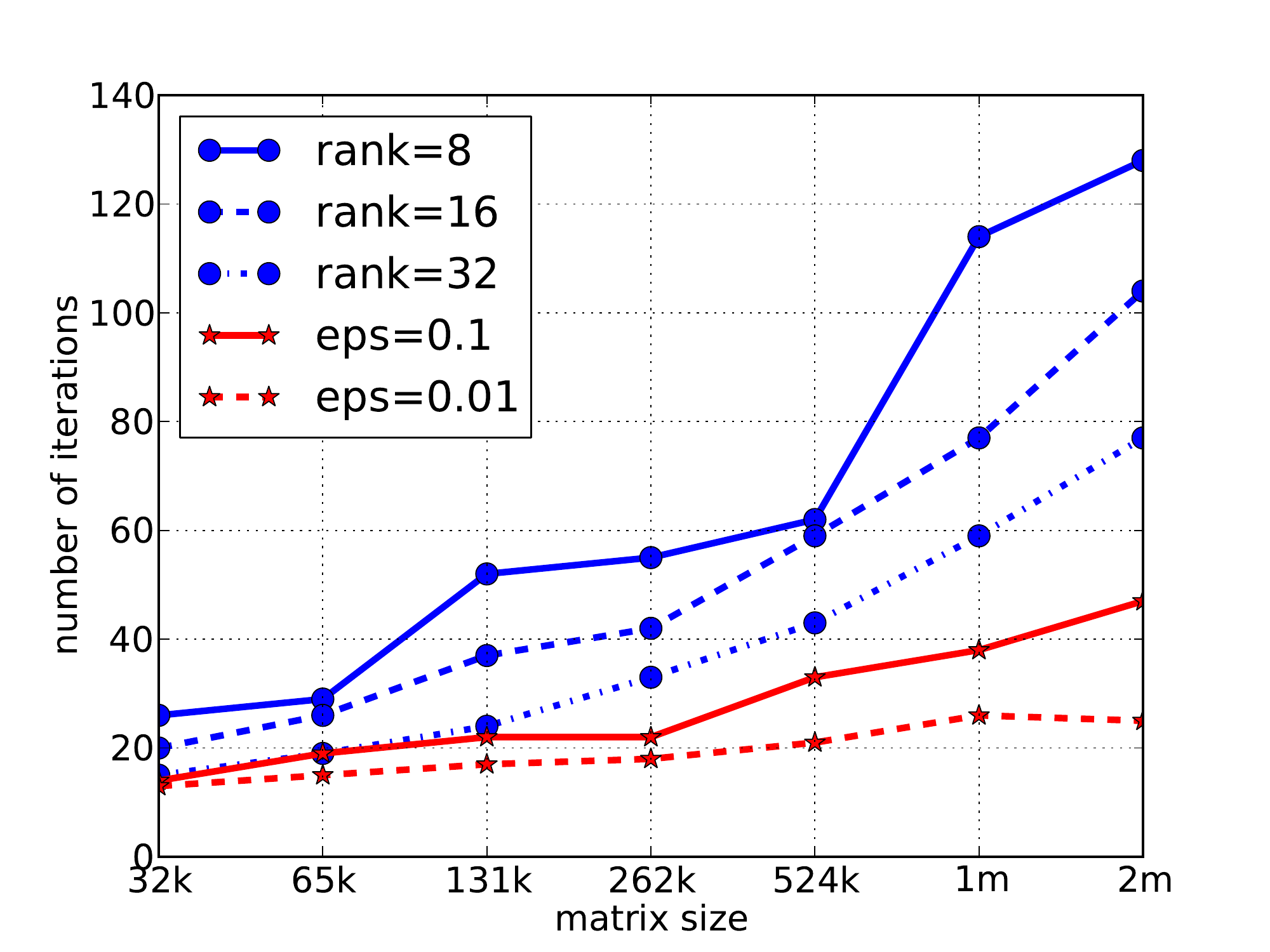}}
  \subfigure[Total time (VC-Poisson)]
  {\includegraphics[width=62mm]{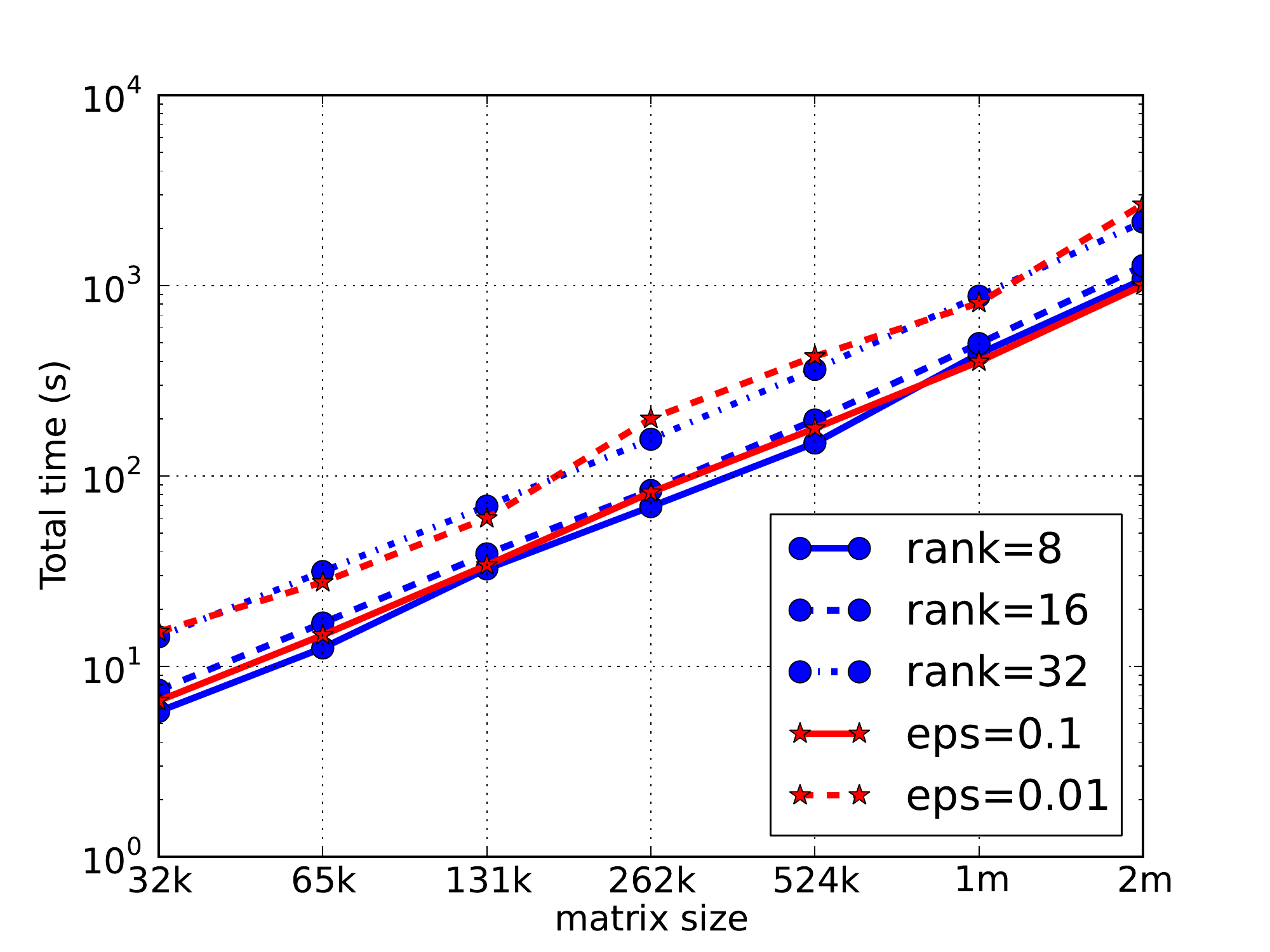}}
  \subfigure[Iteration number (Helmholtz)]
  {\includegraphics[width=62mm]{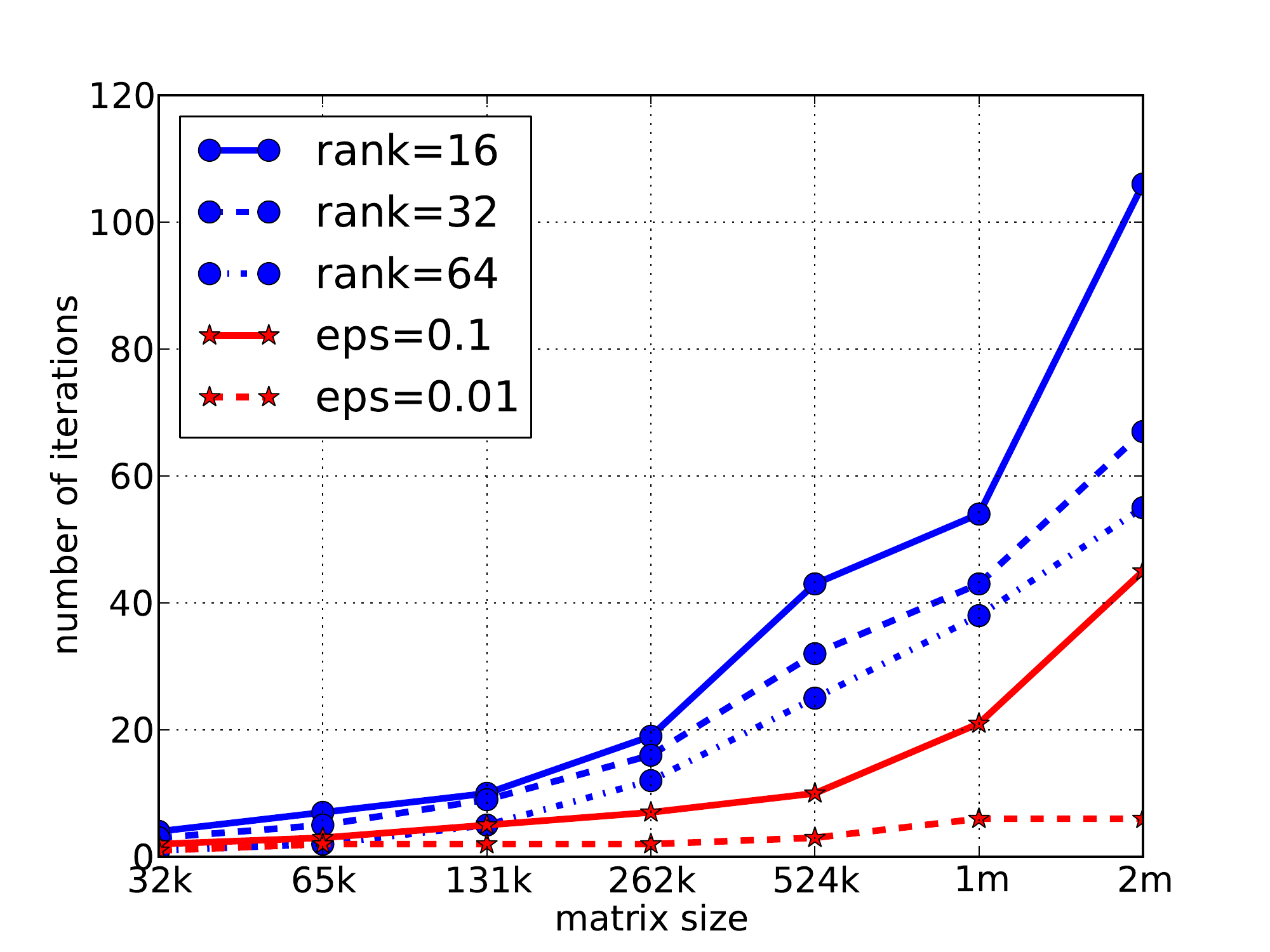}}
  \subfigure[Total time (Helmholtz)]
  {\includegraphics[width=62mm]{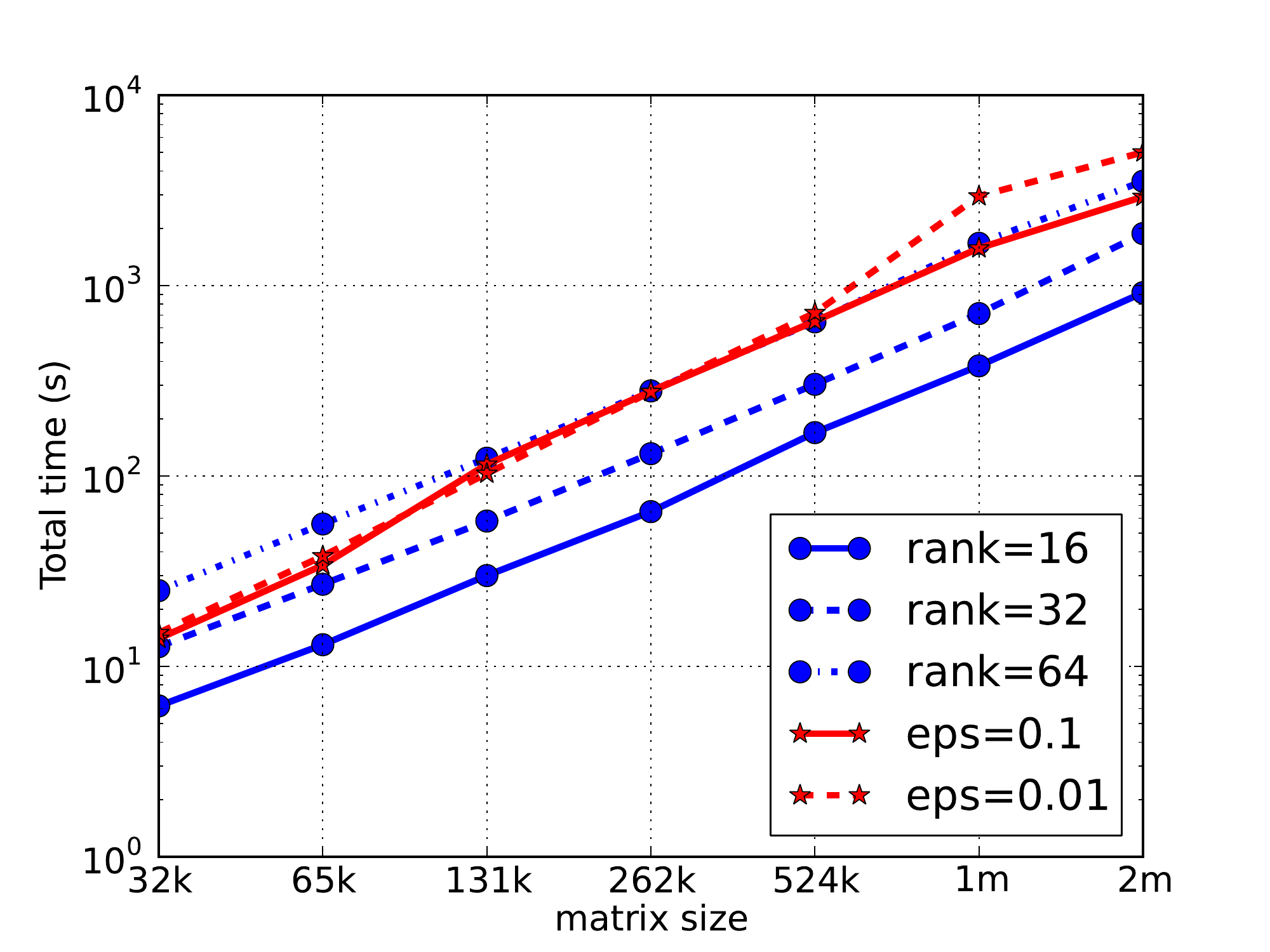}} 
  \caption{Number of iterations and the total time (factorization + solve) for solving three PDEs.  Different low-rank truncation criteria have been used: fixing the rank ${\cal K}$ and fixing the truncation error $\epsilon$. The convergence tolerance is $10^{-12}$ for the first two problems and is $10^{-3}$ for the third problem.}
  \label{fig:iter_P3D}
\end{figure}

We investigated two options to select the rank: (a) dynamic rank based on the singular values of the fill-in blocks, which is computed for a given user-prescribed error tolerance $\epsilon$ (error in the 2-norm), and (b) a fixed value of rank ${\cal K}$. As Fig. \ref{fig:iter_P3D} shows, the former typically gives better quality preconditioners but may be more expensive, while the latter puts a strict upper limit on the memory usage and the factorization cost. In Fig. \ref{fig:iter_P3D}, the total time, which is the sum of the factorization time and the solve time, scales almost linearly with respect to the problem size for all three PDEs.

\subsection{Parallel scalability and scalability bottleneck} \label{subsec:results_parallel}
\begin{figure}[H]
  \centering
  \subfigure[Factorization time (Poisson)]
  {\label{fig:time_P3D}\includegraphics[width=62mm]{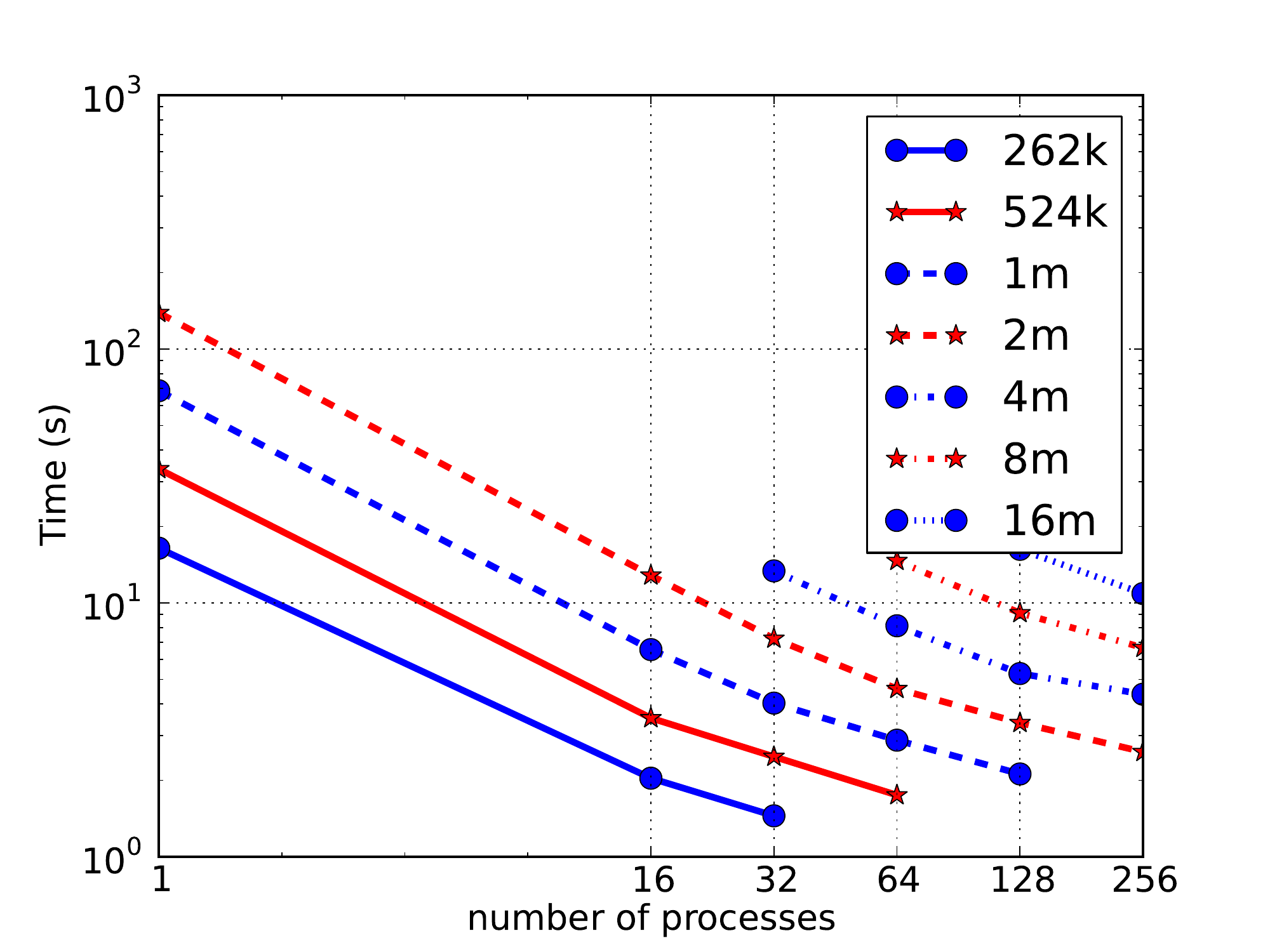}}
  \subfigure[Parallel speedup (Poisson)]
  {\label{fig:speedup_P3D}\includegraphics[width=62mm]{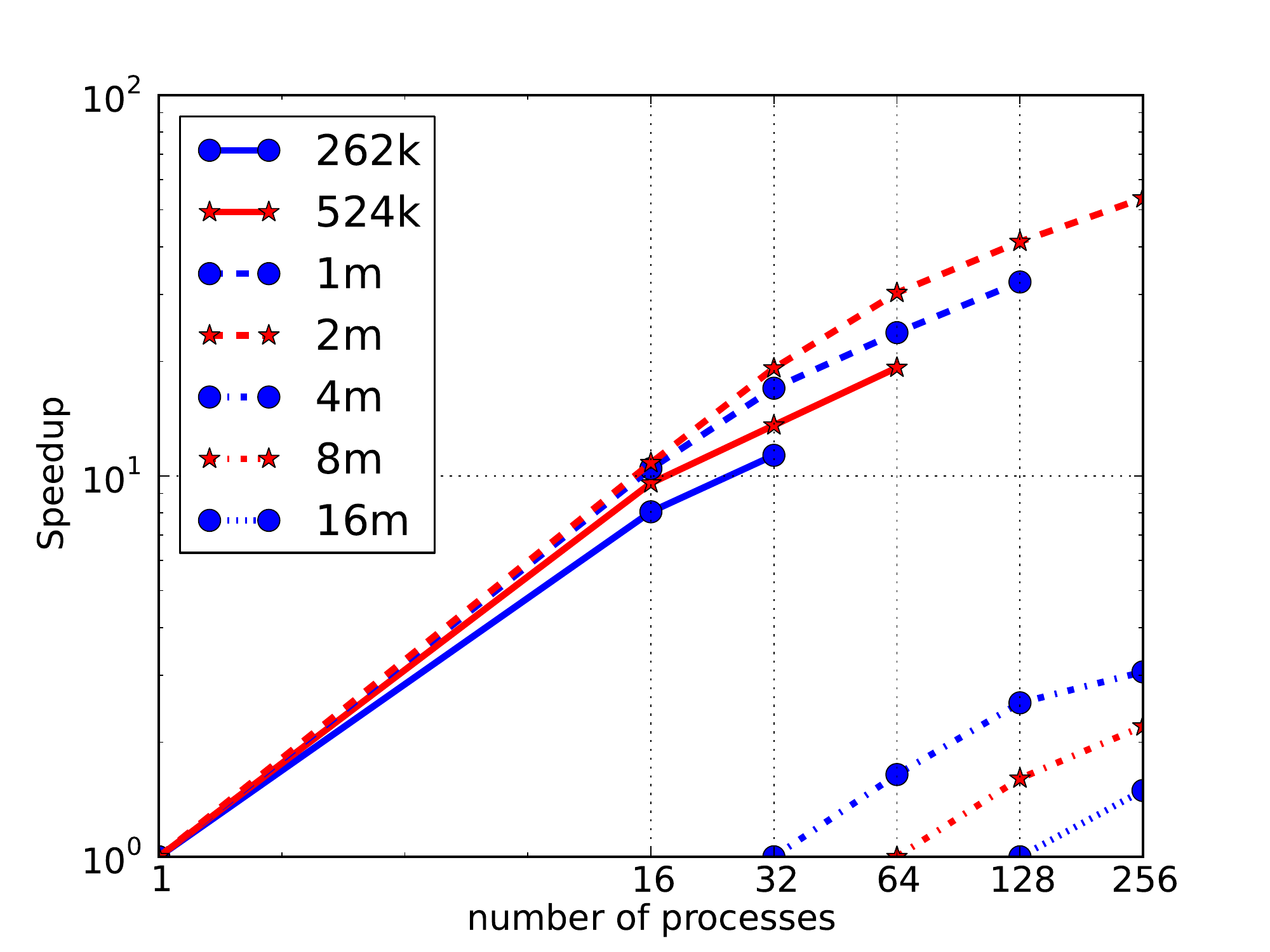}}

  \centering
  \subfigure[Factorization time (VC-Poisson)]
  {\label{fig:time_VP}\includegraphics[width=62mm]{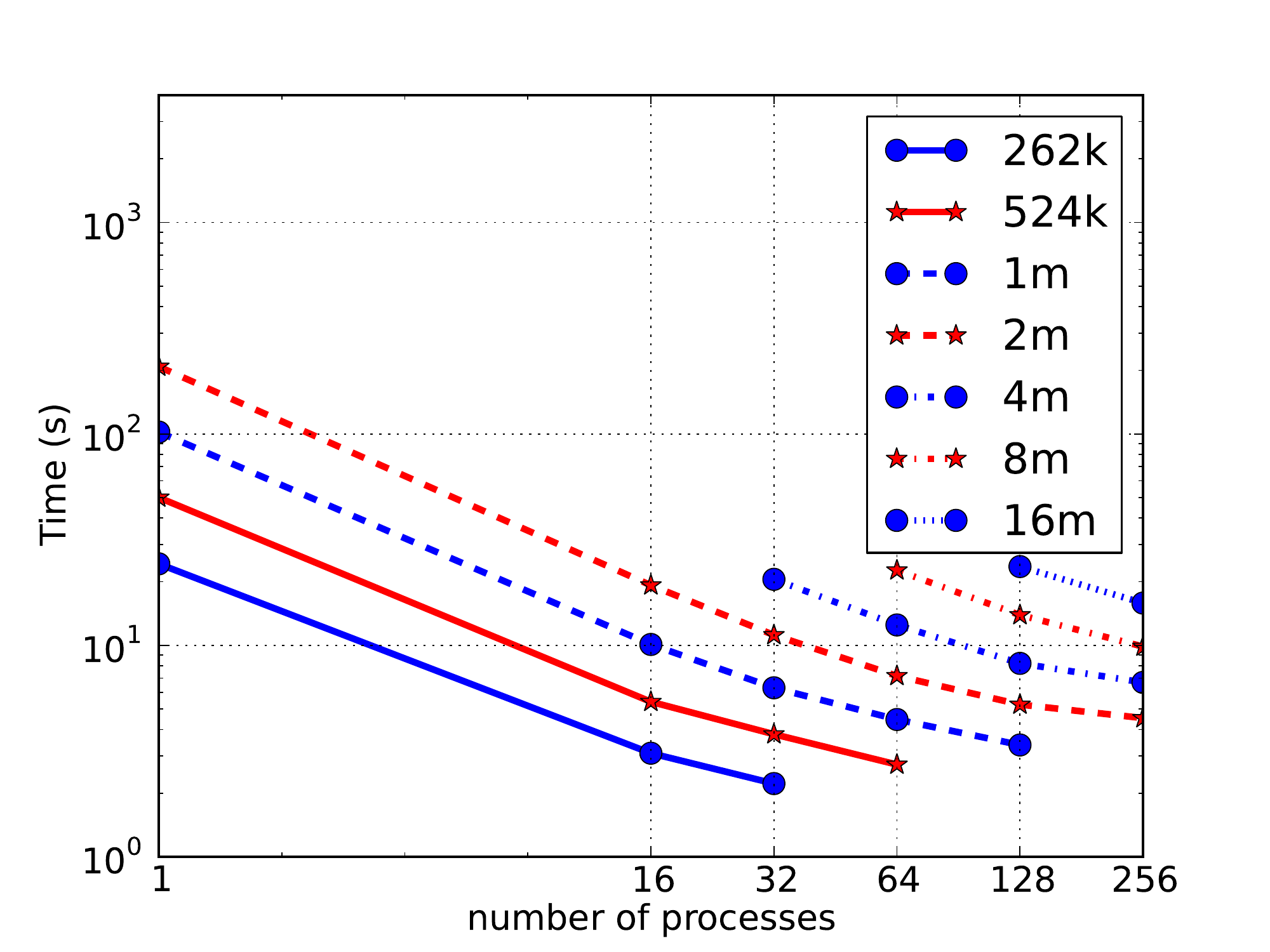}}
  \subfigure[Strong scaling (VC-Poisson)]
  {\label{fig:strong_VP}\includegraphics[width=62mm]{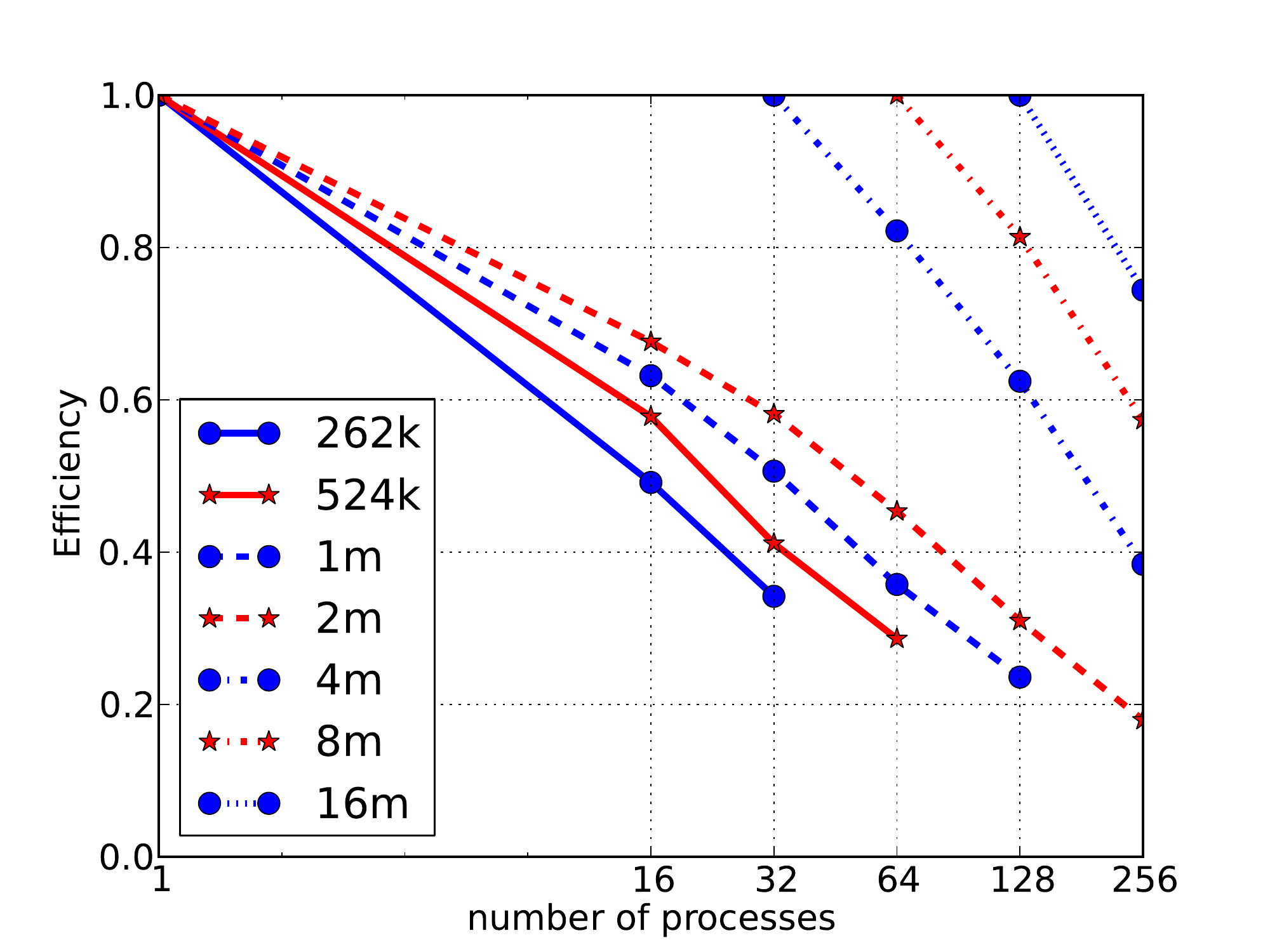}}

  \centering
  \subfigure[Factorization time (Helmholtz)]
  {\label{fig:time_HM}\includegraphics[width=62mm]{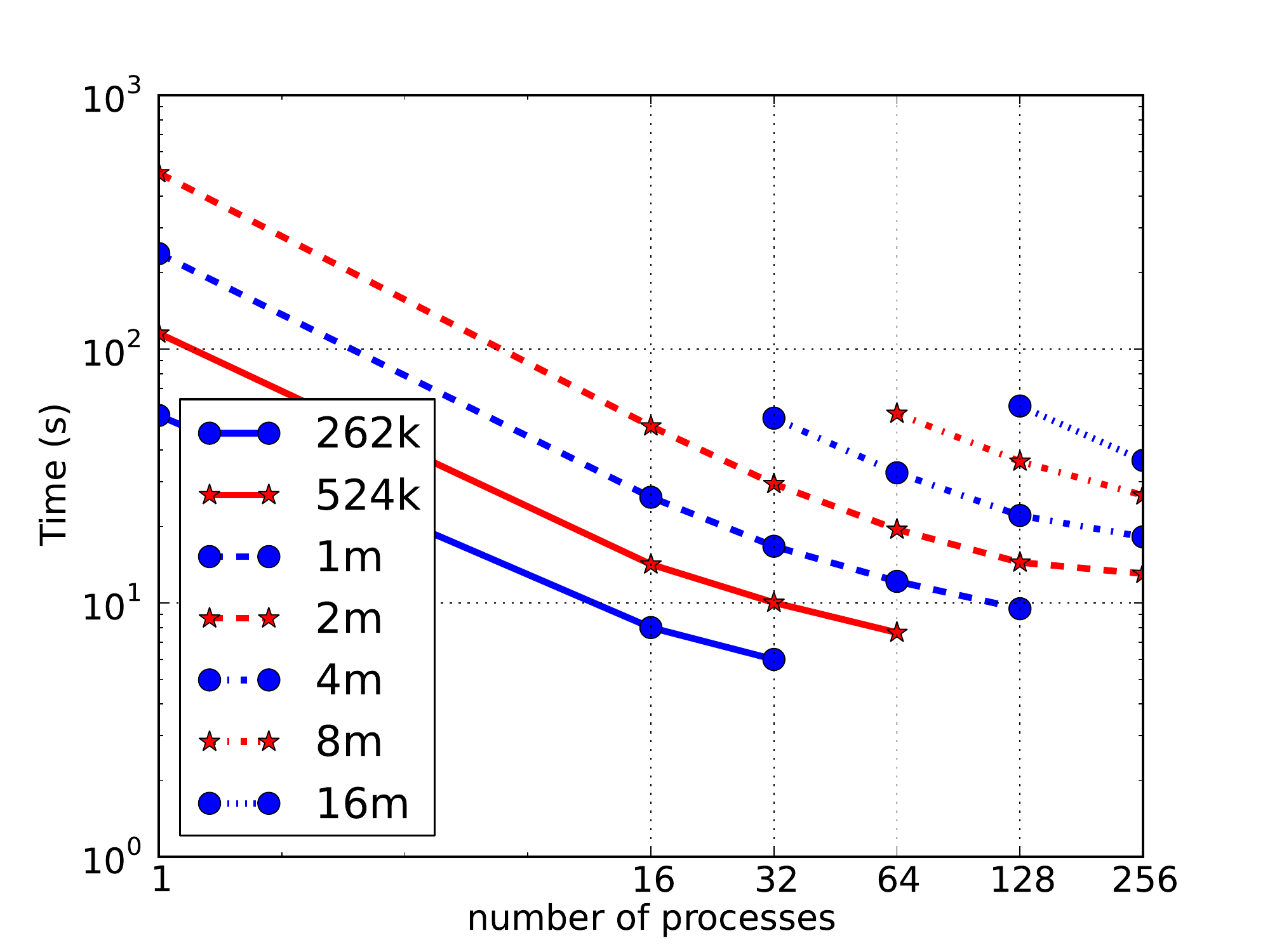}}
  \subfigure[Weak scaling (Helmholtz)]
  {\label{fig:weak_HM}\includegraphics[width=62mm]{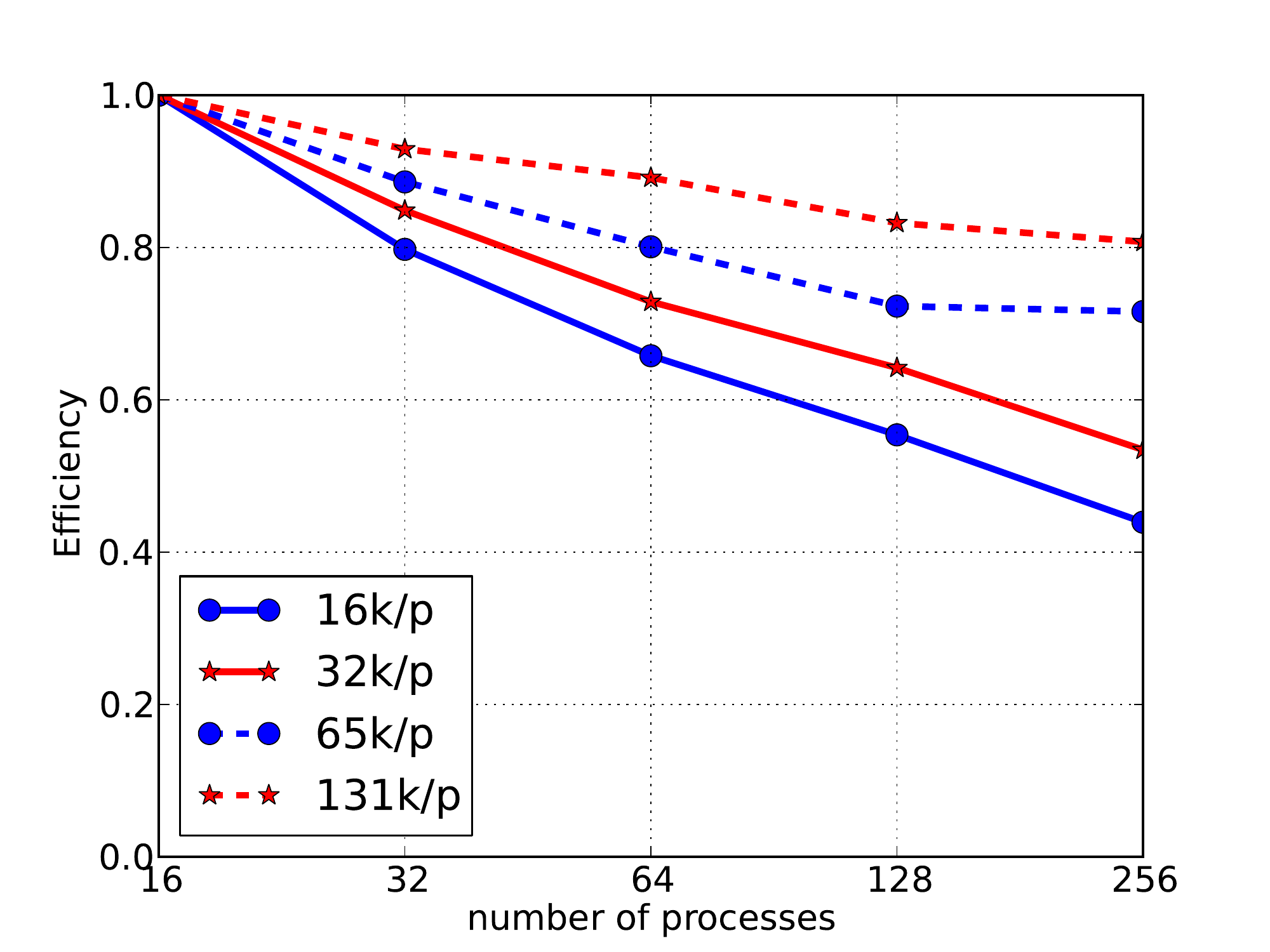}}
  \caption{Parallel scalability results: factorization time and parallel speedups for solving Poisson's equation (${\cal K}=8$), factorization time and strong scaling efficiency for solving the variable-coefficient Poisson's equation (${\cal K}=16$), and factorization time and weak scaling efficiency for solving the Helmholtz equation (${\cal K}=32$).}
  \label{fig:HM}
\end{figure}

In this subsection, we show the parallel scalability of the hierarchical solver. We present both the sequential factorization time and the corresponding parallel results on up to 256 processors (16 processors per node). The definition weak scaling efficiency is given in Section~\ref{sec:analysis}. Fig. \ref{fig:HM} shows the parallel factorization time and the weak scaling efficiency for solving the three PDEs. Note that we stopped the scaling experiments when the number of DOFs per processor is fewer than $8$ thousand. At that scale, the communication costs dominate. The hierarchical solver achieves an average speedup of 45 in factorizing three two-million sized test problems using 256 processors.

\begin{figure}[htbp]
  \centering
  \subfigure[Time breakup for strong scaling]
  {\label{fig:detail_strong_HM}\includegraphics[width=62mm]{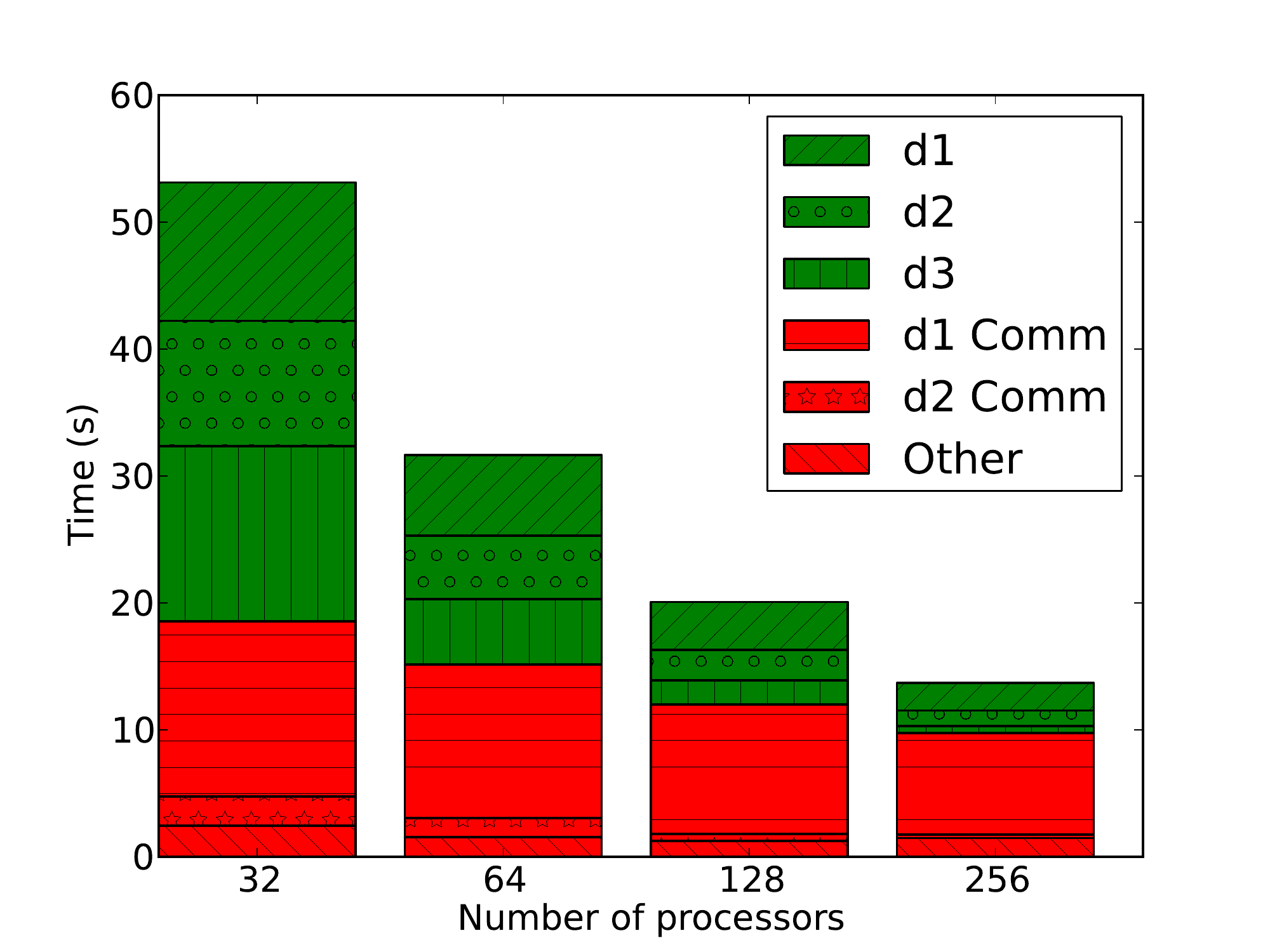}}
  \subfigure[Time breakup for weak scaling]
  {\label{fig:detail_weak_HM}\includegraphics[width=62mm]{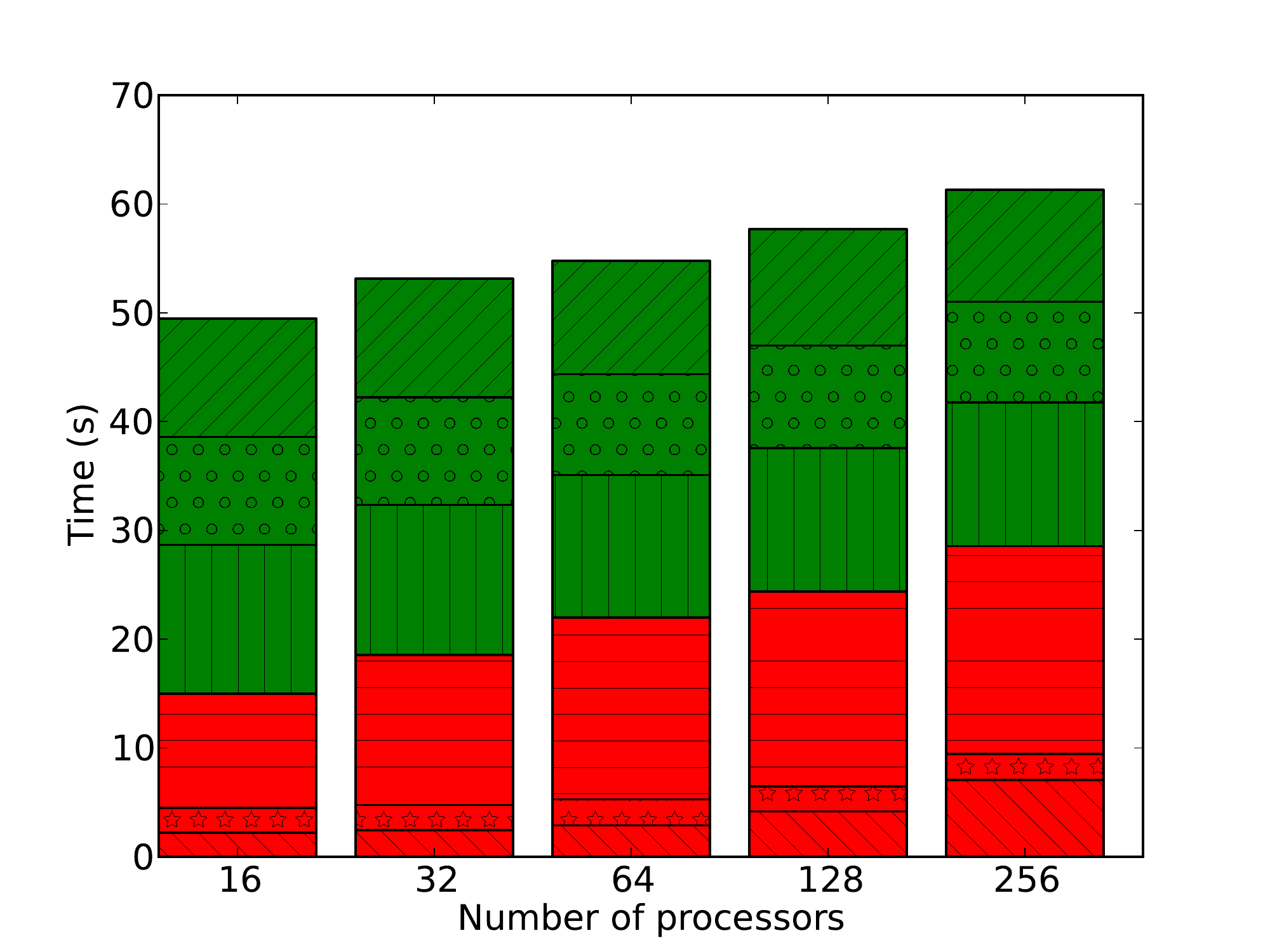}}
  \caption{Breakup of parallel factorization time for solving the Helmholtz equation (${\cal K}=32$). ``d1'', ``d2'' and ``d3'' stand for corresponding computational cost. ``d1 Comm'' and ``d2 Comm'' stand for corresponding communication cost. ``Other'' is mainly the cost of computing d1 coloring.}
  \label{fig:detail_P3D}
\end{figure}

To understand the performance bottleneck, we show breakup of the factorization time in Fig. \ref{fig:detail_P3D}. Since the parallel scalability results are almost the same for solving the three PDEs, we present only the results for solving the Helmholtz equation (${\cal K}=32$). The strong scaling test corresponds to a fixed problem size of four-million DOFs on all processors, and the weak scaling test corresponds to a fixed problem size of 131-thousand DOFs per processor. As Fig. \ref{fig:detail_P3D} shows, in the strong scaling test, the computation time for d1 nodes, d2 nodes and d3 nodes decreases by half when the number of processors doubles, whereas d1 communication time and `Other' time remains almost constant. In the weak scaling test, the computation time for d1 nodes, d2 nodes and d3 nodes remains almost constant when the number of processors increases, whereas d1 communication time and `Other' time increases.

Computing the coloring of d1 nodes affects the parallel scalability in two aspects. First, the time spent in computing the coloring does not scale, as shown in Fig.~\ref{fig:detail_P3D}. Second, more importantly, the quality of the coloring result is poor, which results in more d1 communication than necessary. Recall the hierarchical solver needs a coloring such that every pair of d1 nodes have different colors if they satisfy two conditions: (1) they are within distance two from each other and (2) they belong to two different processors. Unfortunately, we are not aware of any existing algorithm that solves this coloring problem. In our implementation, the standard distance-2 coloring algorithm is used, which ignores the second condition and solves a much more difficult problem. As a result, the output is an unnecessarily stronger coloring result with much more colors than needed. Although the current implementation for computing the coloring is not efficient, the cost can be amortized for solving a sequence of linear systems with the same sparsity pattern, since the coloring result depends only on the graph (symbolic pattern) of the original matrix.

\subsection{Comparison with sparse direct solver} \label{subsec:results_slu}

\begin{figure}[htb]
  \begin{center}
    \subfigure[Running time]
    {{\includegraphics[width=67.5mm]{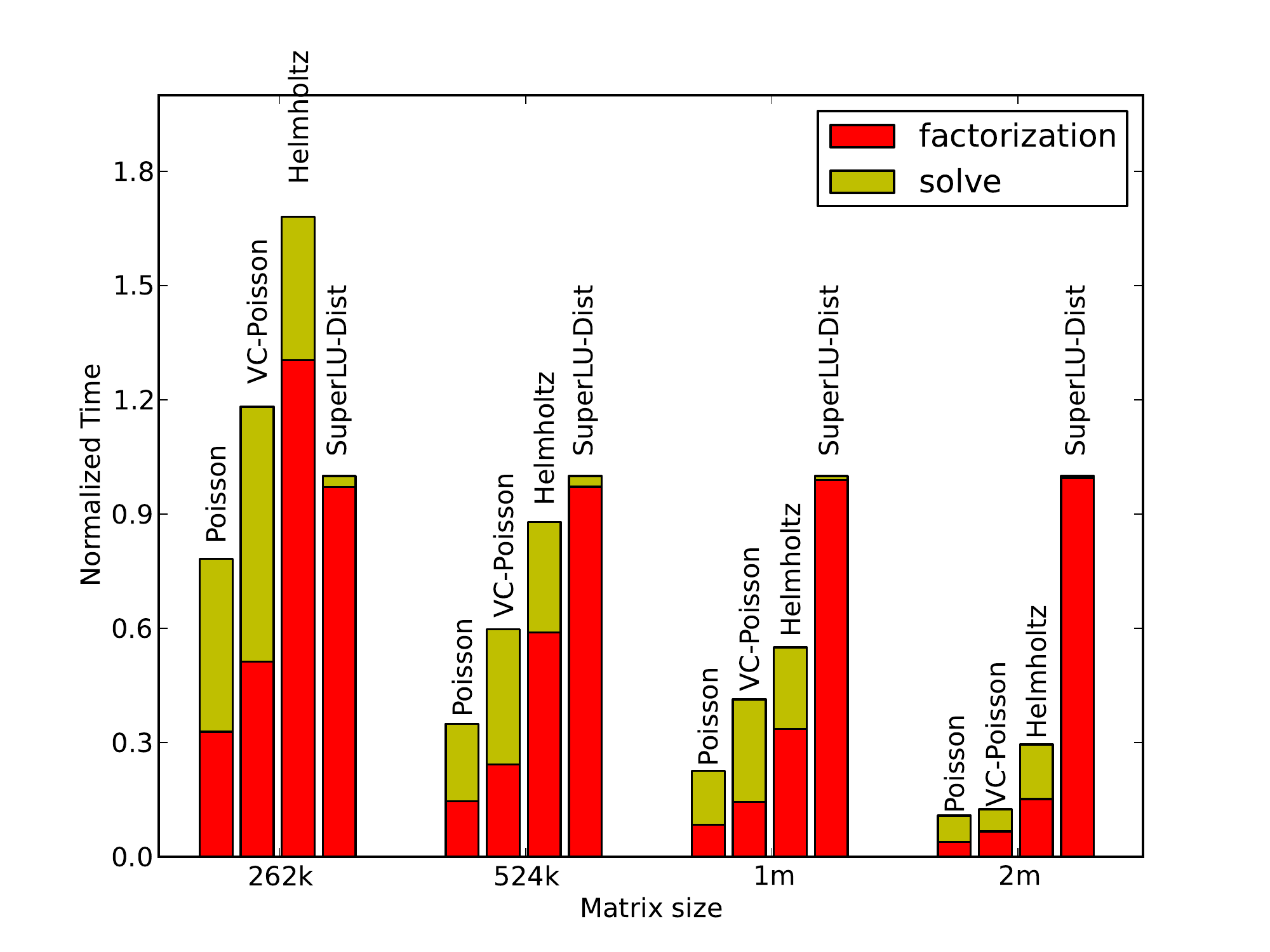}}}
    \subfigure[Memory footprint]
    {{\includegraphics[width=67.5mm]{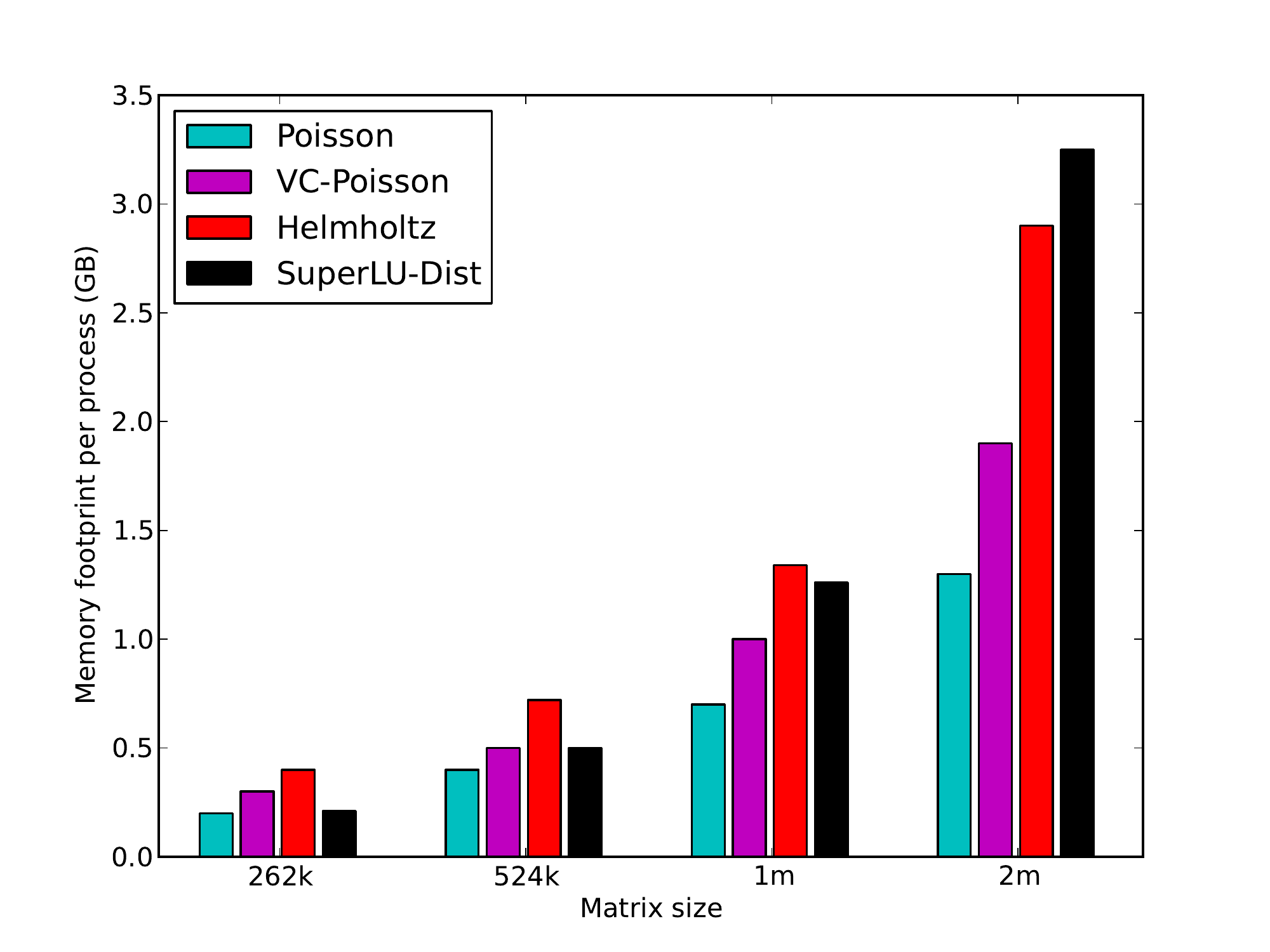}}}
    \caption{Comparison with SuperLU-Dist on 16 processors. For each problem size, timing results are normalized by that of SuperLU-Dist. The timing result and the memory footprint of SuperLU-Dist was almost the same for three PDEs.}
    \label{fig:compare_slu}
  \end{center}
\end{figure}

In this subsection, we present results of comparing our parallel solver to SuperLU-Dist \cite{superlu_ug99,lidemmel03}, a state-of-the-art parallel sparse direct solver, on 16 processors. Timing results and the corresponding memory footprints are shown in Fig. \ref{fig:compare_slu}. For Poisson's equation and the variable-coefficient Poisson's equation, our solver was used as a preconditioner (${\cal K}=8 \textrm{ and } 16$, respectively) for CG with a tolerance of $10^{-12}$. For the Helmholtz equation, our solver was used as a preconditioner (${\cal K}=32$) for GMRES with a tolerance of $10^{-3}$.

In Fig. \ref{fig:compare_slu}, we can see that for a small problem size, e.g., 262k ($128^3$), our solver was not competitive with SuperLU-Dist. But for large problem sizes, our solver was much faster and used less memory than SuperLU-Dist.

Note that the memory footprint of our solver is relatively high. The reason is that in our implementation, (sparse) nonzero blocks in the original sparse linear system are stored using a dense-matrix data structure (for ease of coding), and this choice leads to the large memory consumption. Optimizing this memory usage requires some software-engineering work and is on our to-do list.

\subsection{General linear systems from various applications} \label{subsec:results_ufl}

In this subsection, we present results\footnote{tests in this subsection were run on the Vesper machine (vesper@sandia.gov) at the Sandia National Laboratories. Vesper uses AMD many-core processors and has 128GB of main memory.} for general linear systems from a number of applications, such as electro-physiological model of a torso, numerical weather prediction and atmospheric modeling, geo-mechanical model of earth crust with underground deformation, and gas reservoir simulation for $\text{CO}_2$ sequestration. These problems are from the SuiteSparse Matrix Collection (formerly known as the University of Florida Sparse Matrix Collection) \footnote{https://sparse.tamu.edu/} and the matrices include SPD matrices from unstructured meshes, non-symmetric matrices and indefinite matrices. Important properties of these matrices are shown in \autoref{table:ufl_matrices}.

\begin{table}[!htb]
  \caption{General matrices from various applications}
  \label{table:ufl_matrices}
  \centering
  \begin{tabular}{ | c | c | c | c | c | c |} \hline
    Matrices  &      size& \# nonzero& \thead{symbolic pattern\\ symmetry} & \thead{numeric value\\ symmetry} & \thead{positive\\ definite}  \\ \hline
    torso3    & 0.26M&  4.4M&   no&     no&    no  \\ \hline
    atmosmodd & 1.3M&  8.8M&   yes&    no&    no  \\ \hline
    Geo\_1438 & 1.4M& 60.2M&   yes&   yes&   yes  \\ \hline
    Serena    & 1.4M& 64.1M&   yes&   yes&   yes  \\ \hline
  \end{tabular}
\end{table}

We compared the timing results of SuperLU-Dist and our parallel solver for five low-rank truncation epsilons $(\epsilon = 0.8,0.4,0.2,0.1  \text{ and } 0.05)$. Our solver is used as a preconditioner for CG/GMRES, depending on if the matrix is SPD or not, and the convergence tolerance is $10^{-12}$.
 
\begin{figure}[htbp]
  \centering
  \includegraphics[width=8cm]{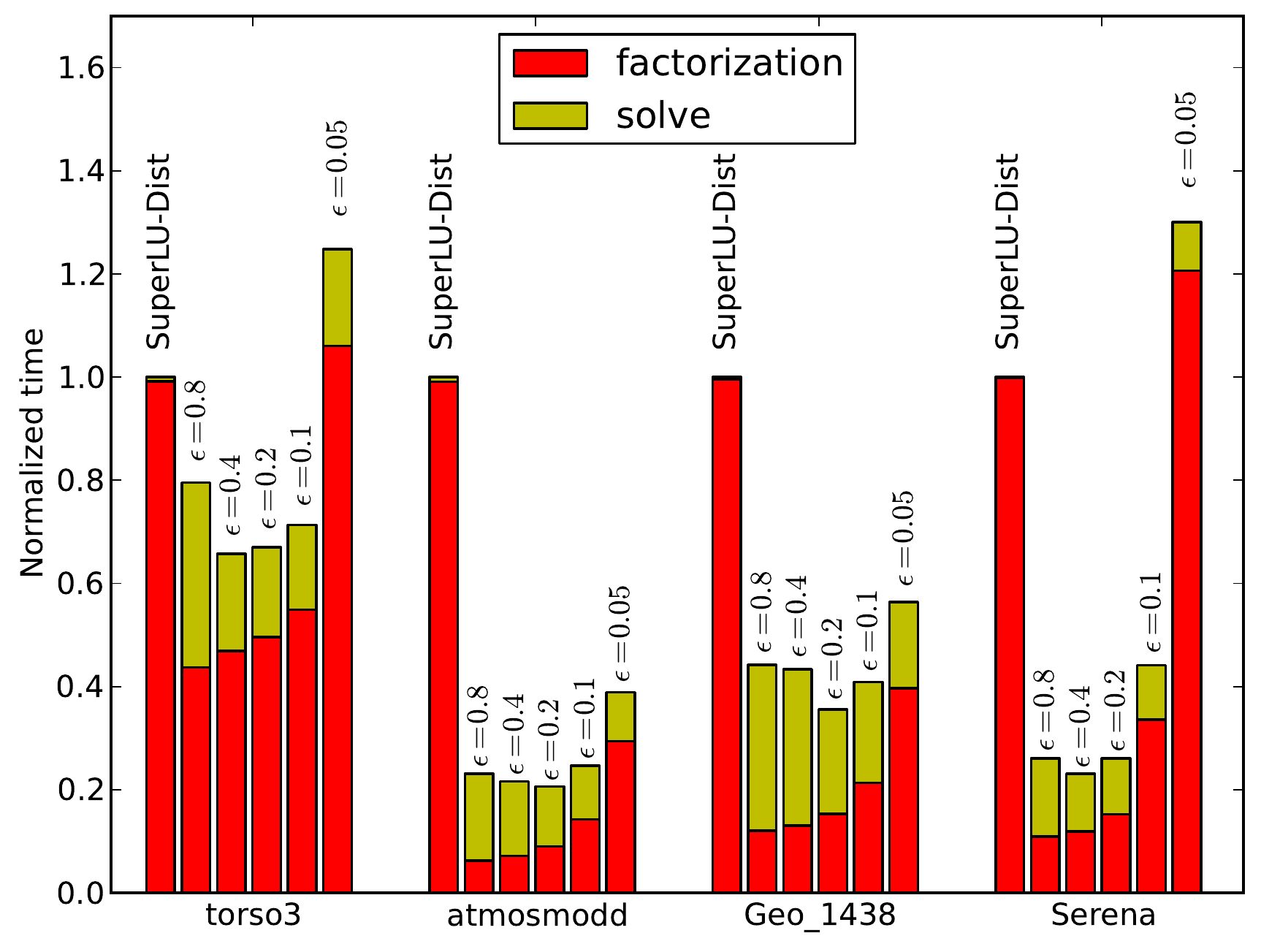}
  \caption{Timing results for different matrices from the University of Florida Sparse Matrix Collection. For each matrix, the results are normalized with respect to the total solve time (factorization + solve) of SuperLU-Dist.}
  \label{fig:ufl_matrices}
\end{figure}

The comparisons are shown in \autoref{fig:ufl_matrices}. In \autoref{fig:ufl_matrices}, as the low-rank truncation error $\epsilon$ decreases (more accurate factorizations), the factorization time increases whereas the solve time typically decreases because of the reduced CG/GMRES iteration numbers. As a result, the total solve time achieves its optimal with $\epsilon$ being around $0.2$ or $0.4$ for these matrices. The speedups of our parallel solver compared with SuperLU-Dist are 1.5x, 4.8x, 2.8x and 4.3x.

\section{Conclusions and future work} \label{sec:conclusion}

We have presented the first parallel algorithm and implementation of the LoRaSp hierarchical solver. Although we derive the algorithm only for SPD matrices in Section~\ref{sec:serial}, the algorithm extends to general matrices. As demonstrated by the numerical experiments, the parallel solver works for non-symmetric matrices, indefinite matrices and matrices from unstructured grids.

The linear complexity conditions are presented in Section~\ref{sec:analysis}. In practice, these conditions may not be satisfied perfectly, but the total solve time (factorization + solve) is observed to scale almost linearly for numerical experiments in Section~\ref{subsec:results_serial}. To solve large linear systems that cannot be stored on a single compute node, distributed memory parallel computing is necessary. We have shown that our solver achieves good speedups on up to 256 processors.

Several directions for future research are as follows.
\begin{itemize}
\item Our graph coloring of the boundary is conservative and other strategies may reduce the number of colors and/or execution time. For example, one could color the processor graph instead of the boundary nodes, which would give more coarse-grained communication but likely more idle time (poor load balance).
\item We plan an MPI+X implementation that exploits thread parallelism for dense linear algebra. We believe such an approach may improve the parallel scalability on a large number of compute nodes.
\item A variant of LoRaSp has been shown to work successfully for dense matrices \cite{coulier2017inverse}. Our parallel algorithm could be extended to solve dense linear systems.
\end{itemize}

\section*{Acknowledgments}
{\small
Chen, Pouransari, and Darve were supported in part by the U.S.~Department of Energy's National Nuclear Security Administration under Award Number DE-NA0002373-1; Boman, Rajamanickam, Chen and Darve were supported in part by Sandia's Laboratory Directed Research and Development (LDRD) program. Sandia National Laboratories is a multimission laboratory managed and operated by National Technology and Engineering Solutions of Sandia, LLC., a wholly owned subsidiary of Honeywell International, Inc., for the U.S. Department of Energy’s National Nuclear Security Administration under contract DE-NA-0003525. This research used resources of the National Energy Research Scientific Computing Center, a DOE Office of Science User Facility supported by the Office of Science of the U.S.~Department of Energy under Contract No. DE-AC02-05CH11231.}





\small
\bibliographystyle{model1-num-names}
\bibliography{ParCompt}







\end{document}